\documentclass[12pt]{article}
\usepackage{amsmath,amssymb}
\usepackage{multirow}
\usepackage{multicol}
\usepackage[usenames,dvipsnames]{xcolor}
\definecolor{light}{rgb}{0.5, 0.5, 0.5}

\usepackage[a4paper, total={6in, 9in}, margin=1in]{geometry}
\usepackage{mathtools}
\usepackage[colorlinks,
            linkcolor=black,
            anchorcolor=black,
            citecolor=black]{hyperref}
\usepackage{natbib}

\DeclarePairedDelimiterX\Basics[1](){ #1}

\usepackage{amsthm}

\usepackage{mathrsfs}
\newtheorem{theorem}{Theorem}[section]

\newtheorem{proposition}[theorem]{Proposition}
\newtheorem{corollary}[theorem]{Corollary}
\usepackage{xcolor}
\usepackage{graphicx}
\graphicspath{ {./images/} }
\usepackage{algorithm} 
\usepackage{algpseudocode} 
\usepackage{caption}
\usepackage{scalerel,stackengine}
\usepackage{blindtext}
\usepackage{authblk}
\usepackage{graphicx}
\newtheorem{condition}{Condition}
\newtheorem{alg}{Algorithm}
\newtheorem{setting}{Setting}
\usepackage{enumitem}
\usepackage{ulem}
\usepackage{booktabs}
\usepackage{multirow}
\usepackage{lscape}
\usepackage{adjustbox}
\usepackage{rotating}
\usepackage{xr}

\def\trans{^{\scriptscriptstyle \sf T}}
\newtheorem{remark}{Remark}

\usepackage[many]{tcolorbox}
\tcbuselibrary{breakable}

\stackMath
\newcommand\reallywidehat[1]{%
\savestack{\tmpbox}{\stretchto{%
  \scaleto{%
    \scalerel*[\widthof{\ensuremath{#1}}]{\kern.1pt\mathchar"0362\kern.1pt}%
    {\rule{0ex}{\textheight}}
  }{\textheight}%
}{2.4ex}}%
\stackon[-6.9pt]{#1}{\tmpbox}%
}
\title{{\LARGE \bf Minimax Rate-Optimal Inference for Individualized Quantile Treatment Effects in High-dimensional Models}}
\date{}
\author{Jiachen Sun and Yin Xia \\
Department of Statistics and Data Science, Fudan University}

\begin{document}

\setlength{\abovedisplayskip}{3pt}
\setlength{\belowdisplayskip}{3pt}

\maketitle{} 

\def\spacingset#1{\renewcommand{\baselinestretch}%
{#1}\small\normalsize} \spacingset{1.2}
\begin{abstract}
The quantification of treatment effects plays an important role in a wide range of applications, including policy making and bio-pharmaceutical research.
In this article, we study the quantile treatment effect (QTE) while addressing two specific types of heterogeneities: (a) personalized heterogeneity, which captures the varying treatment effects for different individuals, and (b) quantile heterogeneity, which accounts for how the impact of covariates varies across different quantile levels. 
A well-designed debiased estimator for the individualized quantile treatment effect (IQTE) is proposed to capture such heterogeneities effectively. 
We show that this estimator converges weakly to a Gaussian process as a function of the quantile levels and propose valid statistical inference methods, including the construction of confidence intervals and the development of hypothesis testing decision rules. 
In addition, the minimax optimality frameworks for these inference procedures are established. Specifically, we derive the minimax optimal rates for the expected length of confidence intervals and the magnitude of the detection boundary for hypothesis testing procedures, illustrating the superiority of the proposed estimator. 
The effectiveness of our methods is demonstrated through extensive simulations and an analysis of the National Health and Nutrition Examination Survey (NHANES) datasets.
\end{abstract}
\noindent\textbf{Keywords}: Debiased estimator; Minimax detection boundary; Minimax expected length; Quantile regression; Weak convergence.

\newpage
\spacingset{1.81}
\section{Introduction}\label{intro.sec}
The evaluation of treatment effects is gaining increasing prominence across a diverse range of sociological and biological fields, such as policy development and bio-pharmaceutical industries {\citep[e.g.][]{imbens2009,kosorok2019}}.
In social science, this includes investigating the causal impact of past experience on current brand preferences \citep{Bart2012} and analyzing the effects of upper secondary schooling policies on hourly wages in Indonesia \citep{SASAKI2023394}. 
Similarly, in bio-pharmaceutical applications, treatment effect evaluation is essential for comparing the effectiveness of non-nucleoside reverse transcriptase inhibitors, such as rilpivirine and efavirenz for HIV-infected patients \citep{liu2017}, as well as assessing the efficacy of remdesivir versus convalescent plasma for COVID-19 treatments \citep{mccaw2020}.

Given the critical importance of treatment assessments, there is a growing body of literature focused on accurately estimating these effects and making valid statistical inference \citep[e.g.,][]{Hirano2003,Abadie2006,kosuke2013}. 
In many applications, a critical phenomenon is that treatment effects may vary across different individuals \citep[e.g.][]{gold2002,sabine2005}. Therefore, a wealth of methods have been developed in various fields to select or design individualized treatments more effectively.
For instance, \cite{cai2011} employs an individual's physiological information to make personalized treatment selection for modern precision medicine;
\cite{wang2018} constructs optimal individualized treatment rules to guide the utilization of first-line insulin therapy based on the specific characteristics of different diabetes patients. 

Nevertheless, the aforementioned works only take into account the individualized heterogeneity through average treatment effect (ATE) estimations and may overlook the quantile-specific heterogeneity that is evident in many scientific fields \citep[e.g.,][]{Abadie2002}.
As an example, a recent study examining environmental factors and body mass index (BMI) finds that while the middle distributions of BMI are similar between breastfed and formula-fed children aged 5-6 years, significant differences are present at the upper and lower quantiles \citep{beyerlein2011}. Several other findings \citep[e.g.][]{briollais2014, wang2022} also suggest that quantile-based approaches are more appropriate and promising for assessing heterogeneous treatment effects.

\subsection{Related Works}
By incorporating both individualized and quantile heterogeneities, this article addresses the estimation and inference of treatment effects through a quantile regression framework.
Specifically, we adopt the following high dimensional quantile regression models respectively for two treatment groups,
$$
F^{-1}_{{Y}_{k,i}\mid\mathbf{X}_{k,i}}(\tau\vert\mathbf{X}_{k,i})=\mathbf{X}_{k,i}\trans \boldsymbol{\beta}_{\tau,k},\text{ for } 1\leq i \leq n_k, \ k=1,2, \text{ and }\tau\in (0,1),
$$
where $F^{-1}_{{Y}_{k,i}\mid\mathbf{X}_{k,i}}(\tau\vert\mathbf{X}_{k,i})$ denotes the conditional $\tau$-quantile function of ${Y}_{k,i}$ given $\mathbf{X}_{k, i}$, ${Y}_{k,i}\in\mathbb{R}$ and $\mathbf{X}_{k,i}\in\mathbb{R}^{p}$ are the response and the covariates collected from the $i$-th individual within group $k$, and $\boldsymbol{\beta}_{\tau,k} \in \mathbb{R}^p$ with $p \gg n_k $ represents the high-dimensional quantile regression coefficient.

Quantile regression is a well-established approach for investigating the impact of covariates on the conditional distribution of responses and serves as an effective tool for characterizing heterogeneous effects. 
In the literature, a rich body of approaches have been proposed for quantile regression inference, 
{including the classical regression rank score method for linear quantile models \citep{guten1992}, the Wald-type approach for partially {linear quantile models} \citep{wang2009}, direct methods that avoid conditional density estimations in both semi-parametric and non-parametric quantile models \citep{fan2016}, and several Bayesian approaches \citep{yang2012, yang2016};} see a thorough survey in \cite{koenker2017}.
Based on the solid theoretical results therein, numerous studies have emerged focusing on inference for the quantile treatment effect (QTE) in conventional low-dimensional settings.
This includes hypothesis testing for the significance of QTE \cite[e.g.][]{koenker2010,Sun2021} and the estimation of QTE under different circumstances \cite[][among many others]{Chernozhukov2005,firpo2007,lanwang2018}.

However, high-dimensionality presents significant challenges for regression analysis, and regularized estimators have been introduced to ensure consistent estimation and prediction \citep[e.g.][]{lasso,candes2007,Bickel2009,peter2011}. 
Beyond linear regression models, regularized methods for high-dimensional quantile regression have also been extensively developed, addressing problems such as parameter estimation and variable selection \citep[e.g.][]{Belloni2011, bradic2011, fan2014,zheng2015}. 
Despite these advances, regularized approaches often introduce non-negligible biases. To address this, debiasing procedures have been developed to achieve valid statistical inference; see \cite{cai2023statistical} for a comprehensive survey and many references therein. Notably, there has been substantial progress in bias correction for high-dimensional quantile regression.  For example, \cite{zhao2015general} proposes a debiased estimator of the quantile regression parameter that modifies the proposal in \cite{JMLR:v15:javanmard14a}. 
\cite{bradic2017uniform} extends this to uniform statistical inference in the framework of quantile processes, and \cite{zhao2020} further expands these results to a distributed setting. 
Additionally, \cite{yan2023} integrates debiasing estimation with convolution-type smoothed quantile regression to enhance computational efficiency.

Leveraging these debiasing approaches, significant advancements have been achieved in high-dimensional treatment effects estimation and inference, as evidenced by \cite{athey2018, wang2020} for the ATE and \cite{zhang2020, chao2022} for the QTE.
Nevertheless, the above studies often overlook the personalized heterogeneity that is prevalent in many applications.
To accommodate such heterogeneity, \cite{cai2021} proposes a debiased estimator for the individualized average treatment effect (IATE) and obtains its minimax optimality; \cite{giessing2023debiased} targets inference on the individualized quantile treatment effect (IQTE) and develops valid procedures based on a regression rank-score debiased estimator.
Despite these advancements, the statistical minimaxity for the IQTE remains unknown, and the existing quantile methods impose strong structural assumptions on covariates, which limits their applicability for studying the IQTE with general individualized predictors.

\subsection{Our Contributions}
In this work, we focus on inference for the IQTE defined as $\Delta_{\tau,\text {new}}=\mathbf{x}_{\text {new}}\trans(\boldsymbol{\beta}_{\tau,1}-\boldsymbol{\beta}_{\tau,2}),$ where $\mathbf{x}_{\text {new}} \in \mathbb{R}^p$ is a structure-free loading vector that characterizes the features of a new individual. 
We propose novel debiasing estimators for linear functionals of quantile regression coefficients $\boldsymbol{\beta}_{\tau,k}$, specifically $\mathbf{x}_{\text {new}}\trans\boldsymbol{\beta}_{\tau,k}$, $k=1,2$, to study the IQTE and facilitate valid inference.
A crucial component of our approach involves constructing a projection direction based on a new variance-enhancement constraint. This constraint provides a lower bound for the variance of the proposed estimator and is instrumental in balancing bias and variance for general loading vectors. 
Based on such projection, we develop a debiased estimator $\widehat{\mathbf{x}_{\text {new}}\trans \boldsymbol{\beta}_{\tau,k}}$ using both a regularized coefficient estimator and a class of initial sparsity function estimators.
We demonstrate that the proposed debiased estimator converges weakly to a Gaussian process over a compact set of quantile levels. 
Leveraging this result, we establish the asymptotic normality for the estimator of $\Delta_{\tau,\text {new}}$ and construct the corresponding confidence intervals and hypothesis testing procedures.  Finally, we derive the minimax optimality theories for both inference procedures, demonstrating the superiority of our proposed debiased estimator.

Our proposal differs from existing solutions and makes several useful contributions.
First, we debias the IQTE $\Delta_{\tau,\text {new}}$ as a whole by incorporating the individualized loading vector to asymptotically eliminate the bias, providing advantages over direct plug-in estimators that focus on debiasing coefficients \citep[e.g.,][]{zhao2015general, bradic2017uniform} rather than linear functionals.
Second, all the existing literatures on the inference of quantile linear functionals assume $\ell_1$-convergent or $\ell_2$-bounded covariates $\mathbf{x}_{\text {new}}$ \citep[e.g.,][]{bradic2017uniform,yan2023,giessing2023debiased}. In contrast, through the construction of a novel variance-enhancement projection direction, the proposed estimator enjoys desirable properties without any structural assumptions on the loading vector $\mathbf{x}_{\text {new}}$.
Third, we establish the uniform convergence of $\widehat{\mathbf{x}_{\text {new}}\trans \boldsymbol{\beta}_{\tau,k}}$ in the high-dimensional setting, where classical Donsker theorems no longer apply, while imposing fewer restrictions on the covariates than existing studies \citep[e.g.,][]{chao2017,chao2019}.
Fourth, we develop a minimax framework for confidence intervals of the IQTE $\Delta_{\tau,\text {new}}$, which, to the best of our knowledge, is not available in the current literature.
Finally, a minimax optimal rate of the detection boundary is developed for the hypothesis testing of $\Delta_{\tau,\text {new}}$. This result generalizes the minimax framework for the IATE \citep{cai2021}, while relaxing their theoretical foundations under normal distributions \citep{Cai_2012}. Instead, we develop a new set of theoretical tools that target at densities in the functional space with bounded derivatives.

\subsection{Outline of the Paper}
The rest of this article is organized as follows. 
Section \ref{method.sec} presents the problem framework and introduces the debiased estimator for the IQTE, along with the construction of confidence intervals and hypothesis testing procedures. 
Section \ref{theory.sec} provides theoretical analyses that ensure the validity of these procedures and develops minimax optimality results for both confidence intervals and hypothesis tests. 
Sections \ref{simu.sec} and \ref{realdata.sec} evaluate the numerical performance of the proposed procedures through simulations and the NHANES datasets.
All technical lemmas, proofs, and additional numerical results are relegated to the Supplementary Material.

\subsection{Notation}
For a matrix $\mathbf{A}=({A}_{i j})\in \mathbb{R}^{n \times p}$, denote by ${A}_{i j}$ its $(i, j){\text{-th}}$ entry. If $\mathbf{A}$ is symmetric, let $\lambda_{\min }(\mathbf{A})$ and $\lambda_{\max }(\mathbf{A})$ represent the smallest and largest eigenvalues of $\mathbf{A}$, respectively. 
For a vector $\boldsymbol{\boldsymbol{\xi}} \in \mathbb{R}^p$, denote by $\operatorname{supp}(\boldsymbol{\xi})$ the support of $\boldsymbol{\xi}$ and let $\|\boldsymbol{\xi}\|_0^{}=\text{Card}\{\text{supp}(\boldsymbol{\xi})\}$ be the $\ell_0$ norm of $\boldsymbol{\xi}$, where Card$\{\cdot\}$ represents the cardinality of a set. Let $\|\boldsymbol{\xi}\|_q^{}=(\sum\limits_{j=1}^p\left|{\xi}_j\right|^q)^{1/q}$ be the $\ell_q$ norm of $\boldsymbol{\xi}$ for $q > 0$ and define the $\ell_{\infty}$ norm of $\boldsymbol{\xi}$ by $\|\boldsymbol{\xi}\|_{\infty}^{}=\max\limits_{1 \leq j \leq p}\left|{\xi}_j\right|$.
Denote by $C,C_0,c_0,C_1,c_1,\ldots$ some universal positive constants that may vary from place to place. 
For two positive sequences $\{a_n\}$ and $\{b_n\}$, write $a_n \lesssim b_n$ if $a_n \leq C b_n$ for all sufficiently large $n$, and $a_n \asymp b_n$ if $a_n \lesssim b_n$ and $b_n \lesssim a_n$; write $a_n \ll b_n$ if $\limsup\limits_{ n\rightarrow \infty} {a_n}/{b_n}=0$. 
Let $a \wedge b=\min\{a,b\}$, for $a,b\in\mathbb{R}$. 
The sub-Gaussian norm of a random variable $T$ is defined as $\|{T}\|_{\psi_2}=\sup\limits_{q \geq 1} \frac{1}{\sqrt{q}}(\mathbb{E}|{T}|^q)^{1/q}$, and the sub-Gaussian norm of a random vector $\mathbf{T} \in \mathbb{R}^p$ is defined as $\|\mathbf{T}\|_{\psi_2}=\sup\limits_{\mathbf{v} \in \mathbb{S}^{p-1}}\|\langle \mathbf{v}, \mathbf{T}\rangle\|_{\psi_2}$, where $\mathbb{S}^{p-1}$ is the unit sphere in $\mathbb{R}^p$.
{Let $\mathbb{I}\left\{\cdot\right\}$ denote the indicator function.}
For any integer $l\geq 0$, let $g^{(l)}(\cdot)$ denote the $l$-th derivative of an $l$-th differentiable function $g(\cdot)$, with $g^{(0)}(\cdot)=g(\cdot)$. 

\section{Methodology}\label{method.sec}

Recall that, we are interested in the following high-dimensional quantile regression models for two treatment groups,
\begin{equation}\label{model}
	F^{-1}_{{Y}_{k,i}\mid\mathbf{X}_{k,i}}(\tau\vert\mathbf{X}_{k,i})=\mathbf{X}_{k,i}\trans \boldsymbol{\beta}_{\tau,k},\text{ for }  1\leq i \leq n_k,\ k=1,2,  \text{ and }\tau\in(0,1),
\end{equation}
where $F^{-1}_{{Y}_{k,i}\mid\mathbf{X}_{k,i}}(\tau\vert\mathbf{X}_{k,i})=\inf\left\{y \in \mathbb{R}: \mathbb{P}\left(Y_{k,i}\leq y\mid \mathbf{X}_{k,i}\right) \geq \tau \right\}$ denotes the conditional $\tau$-quantile function of ${Y}_{k,i}$ given $\mathbf{X}_{k, i}$.
Additionally, $\mathbf{X}_k=\left(\mathbf{X}_{k, 1}, \ldots, \mathbf{X}_{k, n_k}\right)\trans \in \mathbb{R}^{n_k\times p}$ represents independent and identically distributed (i.i.d.) covariates from group $k$, $\mathbf{Y}_k=\left(Y_{k,1}, \ldots, Y_{k, n_k}\right)\trans \in \mathbb{R}^{n_k}$ denotes the corresponding responses, and $\boldsymbol{\beta}_{\tau,k} \in \mathbb{R}^p$ with $p \gg n_k $ is the high-dimensional coefficient of our interest. 
Without loss of generality, we assume that $\mathbb{E}({  \mathbf{X}_k}) = \mathbf{0}$ for $k=1,2$.
We abbreviate $F^{-1}_{{Y}_{k,i}\mid\mathbf{X}_{k,i}}(\cdot\vert\mathbf{X}_{k,i})$ as $F^{-1}_{k,i}(\cdot)$, and denote $F_{k,i}$ and $f_{k,i}$ as the conditional cumulative distribution function (CDF) and density function of the response $Y_{k,i}$ given $\mathbf{X}_{k, i}$, respectively.

Note that, given an individualized covariate $\mathbf{x}_{\text {new}}\in \mathbb{R}^p$, the $\tau$-quantile of the response is $\mathbf{x}_{\text {new}}\trans\boldsymbol{\beta}_{\tau,k}$ for treatment group $k$. 
Hence, to measure the IQTE for a fixed quantile level $\tau\in \mathcal{U} \subset (0,1)$, where $\mathcal{U}$ is a compact set, one is interested in the following quantity:
$$\Delta_{\tau,\text {new}}=\mathbf{x}_{\text {new}}\trans\left(\boldsymbol{\beta}_{\tau,1}-\boldsymbol{\beta}_{\tau,2}\right).$$
From the definition of $\Delta_{\tau,\text {new}}$, we can clearly observe the aforementioned two aspects of heterogeneities: the personalized heterogeneity captured by the inclusion of $\mathbf{x}_{\text {new}}$, and the quantile heterogeneity addressed by selecting different quantile levels.

In this article, for a given covariate $\mathbf{x}_{\text {new}}$, we aim to propose a consistent estimator for $\Delta_{\tau,\text {new}}$ and subsequently construct the corresponding confidence interval. 
In addition, we are interested in studying the following hypothesis testing problem that compares the treatments between the two groups:
$$H_0:\Delta_{\tau,\text {new}}\leq0\quad\quad v.s.\quad\quad H_1:\Delta_{\tau,\text {new}}>0.$$
The detailed estimation and inference procedures will be provided next.

\subsection{Estimation for the IQTE \texorpdfstring{$\Delta_{\tau,\text {new}}$}{DeltaNew}}
As mentioned in Section \ref{intro.sec},  debiasing is a popularly adopted approach to reduce the bias produced by the regularized methods in high-dimensional regression models \citep{zhang2014,vandegeer2014,JMLR:v15:javanmard14a}.
Intuitively, one can {directly} plug in a debiased estimator of $\boldsymbol{\beta}_{\tau,k}$ and obtain an estimator for $\Delta_{\tau,\text {new}}$ as follows: 
\begin{equation}\label{naive.est}
\widehat{\Delta}_{\tau,\text {new}}^{\scriptscriptstyle{\mathtt{db}}}=\mathbf{x}_{\text {new}}\trans\left(\widehat{\boldsymbol{\beta}}_{\tau,1}^{\scriptscriptstyle{\mathtt{db}}}-\widehat{\boldsymbol{\beta}}_{\tau,2}^{\scriptscriptstyle{\mathtt{db}}}\right),
\end{equation}
where $\widehat{\boldsymbol{\beta}}_{\tau,1}^{\scriptscriptstyle{\mathtt{db}}}$ and $\widehat{\boldsymbol{\beta}}_{\tau,2}^{\scriptscriptstyle{\mathtt{db}}}$ are debiased estimators that can be approached by \cite{zhao2015general,bradic2017uniform,zhao2020}, among many others. 
However, the above estimator $\widehat{\Delta}_{\tau,\text {new}}^{\scriptscriptstyle{\mathtt{db}}}$ neglects the possible influence introduced by the loading vector $\mathbf{x}_{\text {new}}$, which could lead to invalid inference unless stringent assumptions, such as $\ell_1$-sparsity \citep{bradic2017uniform}, are imposed on $\mathbf{x}_{\text {new}}$.
The numerical results in Section \ref{simu.sec} support this concern.
Therefore, in order to make valid inference without any structural assumptions on the loading vector  $\mathbf{x}_{\text {new}}$, we propose to debias $\Delta_{\tau,\text {new}}$ as a whole, specifically by developing a debiased estimator that incorporates $\mathbf{x}_{\text {new}}$ when constructing the projection direction.
Notably, a novel variance-enhancement constraint is proposed, which is critical, as demonstrated in Proposition \ref{prop2} and Remarks \ref{remark.structure-free} and \ref{x_new}.

Next, we outline the proposed estimation procedure in Algorithm \ref{alg1}, followed by detailed explanations and intuitive insights in Section \ref{explanations.sec}, focusing on two key steps: construction of the projection vector and development of the debiasing procedure.

\vspace{0.2in}
\begin{tcolorbox}[breakable, colframe=white, colback=white, enhanced, frame hidden, sharp corners,
    overlay={\draw[line width=0.5pt] (frame.south west) -- (frame.south east);}, 
    overlay first={\draw[line width=0.8pt] (frame.north west) -- (frame.north east);}, 
    left=-1.5mm,         
    right=-1.5mm,        
    coltitle=black, title= \begin{alg} \label{alg1} \textup{A debiased estimator of} $\Delta_{\tau,\textup{new}}$ \end{alg}]
    \vspace{-3mm}
    \noindent \rule[0.5mm]{\textwidth}{0.4pt}
	\textbf{Input:} $\mathbf{x}_{\text {new}}\in \mathbb{R}^p$, $(\mathbf{Y}_k,{  \mathbf{X}_k})$, an initial coefficient estimator $\widehat{\boldsymbol{\beta}}_{\tau,k}$, sparsity function estimators $\widehat{\eta}_{k,i}(\tau)$, and tuning parameters $\lambda_k,\gamma_k$, $1\leq i\leq n_k$, $k=1,2$. 

    \quad\  Step 1. \textbf{Projection construction:} Obtain a projection direction $\widehat{\mathbf{M}}_{\tau,k}$ by
	\begin{align}
		&\widehat{\mathbf{M}}_{\tau,k}=\underset{\mathbf{M}_{\tau,k} \in \mathbb{R}^p}{\arg \min\ } \mathbf{M}_{\tau,k}\trans \widehat{\boldsymbol{\Sigma}}_{\tau,k} \mathbf{M}_{\tau,k}, \label{obj}\\
		\text { s.t. } &\Vert\widehat{\boldsymbol{\Sigma}}_k \mathbf{M}_{\tau,k}-\mathbf{x}_{\text {new}}\Vert_{\infty}^{} \leq\left\|\mathbf{x}_{\text {new}}\right\|_2 \lambda_k, \label{1}\\
        & |\mathbf{x}_{\text {new}}\trans \widehat{\boldsymbol{\Sigma}}_k \mathbf{M}_{\tau,k}-\left\|\mathbf{x}_{\text {new}}\right\|_2^2| \leq\left\|\mathbf{x}_{\text {new}}\right\|_2^2 \lambda_k,\label{2}\\
        &\underset{1\leq i\leq n_k}{\max} |\mathbf{X}_{k, i}\trans\mathbf{M}_{\tau,k} | \leq \left\|\mathbf{x}_{\text {new}}\right\|_2\gamma_k,\label{3}
	\end{align}
	\quad\quad \quad\quad \ \ \  where $\widehat{\boldsymbol{\Sigma}}_k=n_k^{-1} \sum\limits_{i=1}^{n_k} \mathbf{X}_{k, i} \mathbf{X}_{k, i}\trans$ is the sample covariance matrix, and $\widehat{\boldsymbol{\Sigma}}_{\tau,k}=$ 
    
    \quad\quad \quad\quad \ \ \   $n_k^{-1}\sum\limits_{i=1}^{n_k} \left(\widehat{\eta}_{k,i}(\tau)\right)^2\mathbf{X}_{k, i} \mathbf{X}_{k, i}\trans$ is a quantile-adjusted sample covariance matrix.
    \vspace{2mm}

	\quad\  Step 2.\textbf{ Debiasing:} Compute
	\begin{align}
        \quad\quad\quad\quad\widehat{\mathbf{x}_{\text {new}}\trans \boldsymbol{\beta}_{\tau,k}}=\mathbf{x}_{\text {new}}\trans \widehat{\boldsymbol{\beta}}_{\tau,k}+\frac{1}{n_k}\widehat{\mathbf{M}}_{\tau,k}\trans\sum\limits_{i=1}^{n_k} \widehat{\eta}_{k,i}(\tau)\mathbf{X}_{k, i}{\varphi}_{\tau}\left(Y_{k,i}-\mathbf{X}_{k, i}\trans\widehat{\boldsymbol{\beta}}_{\tau,k}\right),\label{4} 
	\end{align} 
	\quad\quad \quad\quad \ \ \ where ${\varphi}_{\tau}(x)=\tau-\mathbb{I}\left\{x\leq0\right\}$ is the quantile score function.
	
	\textbf{Output:} A debiased estimator of the IQTE $\Delta_{\tau,\text {new}}$:
	\begin{align}
		 \widehat{\Delta}_{\tau,\text {new}}=\widehat{\mathbf{x}_{\text {new}}\trans \boldsymbol{\beta}_{\tau,1}}-\widehat{\mathbf{x}_{\text {new}}\trans \boldsymbol{\beta}_{\tau,2}}.\label{deb}
	\end{align}
\end{tcolorbox}

\subsection{Detailed Explanations for Algorithm \ref{alg1}}\label{explanations.sec}

In Step 1, we construct the projection direction $\widehat{\mathbf{M}}_{\tau,k}$ by incorporating one objective function and three functionally distinct constraints, each serving a unique role in validating the proposed estimator.
Specifically, Constraints \eqref{1}  and \eqref{3} are employed to control the biases, while Constraint  \eqref{2} plays a key role in balancing the variance.

More precisely, the estimation error of $\widehat{\mathbf{x}_{\text {new}}\trans \boldsymbol{\beta}_{\tau,k}}$ in \eqref{4} can be decomposed into two main components:
\begin{align}
    \widehat{\mathbf{x}_{\text {new}}\trans \boldsymbol{\beta}_{\tau,k}}-\mathbf{x}_{\text {new}}\trans \boldsymbol{\beta}_{\tau,k}=U_{\tau,k}+\Delta_{\tau,k},\label{dec}
\end{align}
where $U_{\tau,k}$ is a dominating variance term defined by
\begin{align}
U_{\tau,k}=\frac{1}{n_k}\widehat{\mathbf{M}}_{\tau,k}\trans\sum\limits_{i=1}^{n_k}{\frac{1}{f_{k,i}(F_{k,i}^{-1}(\tau))}}\mathbf{X}_{k, i}\varphi_\tau(Y_{k,i}-\mathbf{X}_{k, i}\trans{\boldsymbol{\beta}}_{\tau,k}), \label{var}
\end{align}
and $\Delta_{\tau,k}$ is a bias term defined by
\begin{align}
    \Delta_{\tau,k}=\mathbf{x}_{\text {new}}\trans&(\widehat{\boldsymbol{\beta}}_{\tau,k}-\boldsymbol{\beta}_{\tau,k})+\frac{1}{n_k}\widehat{\mathbf{M}}_{\tau,k}\trans\sum\limits_{i=1}^{n_k}{\widehat{\eta}_{k,i}(\tau)}\mathbf{X}_{k, i}\nu_{\tau,k,i}(\widehat{\boldsymbol{\beta}}_{\tau,k})\nonumber\\
    &+\frac{1}{n_k}\widehat{\mathbf{M}}_{\tau,k}\trans\sum\limits_{i=1}^{n_k}\left(\widehat{\eta}_{k,i}(\tau)-{1}/{f_{k,i}(F_{k,i}^{-1}(\tau))}\right)\mathbf{X}_{k, i}\varphi_\tau(Y_{k,i}-\mathbf{X}_{k, i}\trans{\boldsymbol{\beta}}_{\tau,k}),\label{bias}
\end{align}
with the quantity $\nu_{\tau,k,i}(\widehat{\boldsymbol{\beta}}_{\tau,k})=\mathbb{I}\left\{Y_{k,i}\leq\mathbf{X}_{k, i}\trans{\boldsymbol{\beta}}_{\tau,k}\right\}-\mathbb{I}\{Y_{k,i}\leq\mathbf{X}_{k, i}\trans\widehat{\boldsymbol{\beta}}_{\tau,k}\}$. The initial estimators $\widehat{\boldsymbol{\beta}}_{\tau,k}$ and $\widehat{\eta}_{k,i}(\tau)$ will be specified later in Step 2.

In order to control the bias term $\Delta_{\tau,k}$, we further decompose it into the following {five} refined components:
\begin{align}\label{bias.decomp}
    \Delta_{\tau,k,1}=-(\widehat{\boldsymbol{\Sigma}}_k\widehat{\mathbf{M}}_{\tau,k}-\mathbf{x}_{\text {new}})\trans&{\widehat{\boldsymbol{\zeta}}_{\tau,k}},\Delta_{\tau,k,2}=-\frac{1}{n_k}\widehat{\mathbf{M}}_{\tau,k}\trans\sum\limits_{i=1}^{n_k}(\frac{\widehat{\eta}_{k,i}(\tau)}{1/f_{k,i}(F^{-1}_{k,i}(\tau))}-1)\mathbf{X}_{k, i}\mathbf{X}_{k, i}\trans\widehat{\boldsymbol{\zeta}}_{\tau,k},\nonumber\\
    \Delta_{\tau,k,3}=\frac{1}{n_k}\widehat{\mathbf{M}}_{\tau,k}\trans&\sum\limits_{i=1}^{n_k}\widehat{\eta}_{k,i}(\tau)\mathbf{X}_{k, i}(\nu_{\tau,k,i}(\widehat{\boldsymbol{\beta}}_{\tau,k})-\rho_{\tau,k,i}(\widehat{\boldsymbol{\beta}}_{\tau,k})),\nonumber\\
    \Delta_{\tau,k,4}=-\frac{1}{2n_k}&\widehat{\mathbf{M}}_{\tau,k}\trans\sum\limits_{i=1}^{n_k}{\widehat{\eta}_{k,i}(\tau)}{f^{(1)}_{k,i}(a_{\tau,k,i})}\mathbf{X}_{k, i}(\mathbf{X}_{k, i}\trans\widehat{\boldsymbol{\zeta}}_{\tau,k})^2,\nonumber\\
    {\Delta_{\tau,k,5}=\frac{1}{n_k}\widehat{\mathbf{M}}_{\tau,k}\trans\sum\limits_{i=1}^{n_k}}&{ {\left(\widehat{\eta}_{k,i}(\tau)-{1}/{f_{k,i}(F_{k,i}^{-1}(\tau))}\right)}\mathbf{X}_{k, i}\varphi_\tau(Y_{k,i}-\mathbf{X}_{k, i}\trans{\boldsymbol{\beta}}_{\tau,k}),}
\end{align}
where ${\widehat{\boldsymbol{\zeta}}_{\tau,k}}=\widehat{\boldsymbol{\beta}}_{\tau,k}-\boldsymbol{\beta}_{\tau,k}$, ${\rho_{\tau,k,i}(\widehat{\boldsymbol{\beta}}_{\tau,k})}=F_{k,i}(\mathbf{X}_{k, i}\trans{\boldsymbol{\beta}}_{\tau,k})-F_{k,i}(\mathbf{X}_{k, i}\trans\widehat{\boldsymbol{\beta}}_{\tau,k})$, and $a_{\tau,k,i}$ locates between $\mathbf{X}_{k, i}\trans{\boldsymbol{\beta}}_{\tau,k}$ and $\mathbf{X}_{k, i}\trans\widehat{\boldsymbol{\beta}}_{\tau,k}$. 
In line with the debiasing procedures in high-dimensional statistical inference under both linear regression and quantile regression frameworks \citep[e.g.,][]{cai2021,zhao2015general,bradic2017uniform,giessing2023debiased},
Constraint \eqref{1} is developed to bound the first bias term $\Delta_{\tau,k,1}$, and Constraint \eqref{3} is employed to control the rest four bias terms.

To control the variance term $U_{\tau,k}$, in addition to the objective function that minimizes the quadratic term $\mathbf{M}_{\tau,k}\trans \widehat{\boldsymbol{\Sigma}}_{\tau,k} \mathbf{M}_{\tau,k}$ in order to obtain an upper bound of the variance, we propose a novel Constraint \eqref{2} that additionally provides a lower bound for the variance of the proposed estimator. It is worthy noting that, such lower bound is essential to the development of a unified inference that accommodates both sparse and dense loading vectors $\mathbf{x}_{\text {new}}$. 
Intuitively, without such lower bound, if $\Vert\mathbf{x}_{\text {new}}\Vert_2^{}$ is large, the feasible set of the optimization problem \eqref{obj} might also be large, potentially driving $\widehat{\mathbf{M}}\trans_{\tau,k}\widehat{\boldsymbol{\Sigma}}_{\tau,k}\widehat{\mathbf{M}}_{\tau,k}$ toward zero. This may lead to invalid inference procedures because the variance of the estimators can no longer dominate the bias; see for example, \cite{zhao2015general,bradic2017uniform,giessing2023debiased}, who assume the $\ell_1$-convergent or $\ell_2$-bounded loading vector to guarantee the validity of their procedures.
Henceforth, the proposed Constraint \eqref{2} is crucial for the inference procedures, as it allows for structure-free loading vectors.

Next, in Step 2, based on the projection vector, we perform debiasing through Equation \eqref{4}, which involves two types of initial estimators: $\widehat{\boldsymbol{\beta}}_{\tau,k}$ and $\{\widehat{\eta}_{k,i}(\tau), 1\leq i\leq n_k\}$.
Similarly as the debiasing procedures in linear regression models, $\widehat{\boldsymbol{\beta}}_{\tau,k}$ serves as an initial regularized coefficient estimator and can be approached by many shrinkage methods.
It will be shown in Section \ref{asymp.subsec} that any initial estimator $\widehat{\boldsymbol{\beta}}_{\tau,k}$ satisfying Condition \ref{cond3} can reach the desired theoretical properties and hence can be employed in the construction of $\widehat{\mathbf{x}_{\text {new}}\trans \boldsymbol{\beta}_{\tau,k}}$.
As an example, we can use the approach from \cite{Belloni2011} and estimate $\boldsymbol{\beta}_{\tau,k}$ by fitting the $\ell_1$-penalized quantile regression as 
\begin{align}\label{Belloni.est}
\widehat{\boldsymbol{\beta}}_{\tau,k} \in \underset{\boldsymbol{\beta}_{\tau,k} \in \mathbb{R}^p}{\operatorname{argmin}}\left\{\frac{1}{n_k} \sum_{i=1}^{n_k} {  \varPsi}_\tau\left(Y_{k,i}-\mathbf{X}_{k,i}\trans\boldsymbol{\beta}_{\tau,k} \right)+\lambda\|\boldsymbol{\beta}_{\tau,k} \|_1^{}\right\},
\end{align}
where ${  \varPsi}_\tau(x)=x{  \varphi}_\tau(x)$ is the check function and the regularization parameter satisfies $\lambda\asymp \sqrt{\operatorname{log}p/n_k}$. 
The second class of initial estimators $\widehat{\eta}_{k,i}(\tau)$ serves to estimate the sparsity function $1/f_{k,i}(F_{k,i}^{-1}(\tau))$, which reflects the density of observations at a specific quantile: a smaller value of the sparsity function indicates a higher concentration of the data at this quantile.
Therefore, our target linear functionals at the quantile of interest are easier to be estimated if the corresponding sparsity function is small.
This coincides with the fact that a small sparsity function implies that debiasing is relatively unnecessary, as indicated in Equation \eqref{4}.
Various types of initial estimators can be employed for this purpose, including Koenker-type estimators \citep{koenker_2005,belloni2019}, rank-based estimators \citep{bradic2017uniform}, and kernel-density estimators \citep{tan2022}. 
Depending on the imposed regularity conditions, these estimators achieve consistency at distinct convergence rates. In this work, we adopt the classical Koenker’s quotient estimator \citep{koenker_2005}.
Specifically, based on the relation that
${dF^{-1}_{k,i}(x)}/{dx}={1}/{f_{k,i}(F_{k,i}^{-1}(x))} \text{, for any }x\in\left(0,1\right),$
the sparsity function is estimated by
\begin{align}
    \widehat{\eta}_{k,i}(\tau)=\left\{ {\widehat{F}^{-1}_{k,i}(\tau+h_k)-\widehat{F}^{-1}_{k,i}(\tau-h_k)} \right\}/{2h_k},\label{koen}
\end{align}
where $\widehat{F}^{-1}_{k,i}(\tau+h_k)=\mathbf{X}_{k,i}\trans\widehat{\boldsymbol{\beta}}_{\tau+h_k,k}$ and $\widehat{F}^{-1}_{k,i}(\tau-h_k)=\mathbf{X}_{k,i}\trans\widehat{\boldsymbol{\beta}}_{\tau-h_k,k}$ with $h_k\asymp( (s_k\log p)^2/ n_k)^{1 / 6}$. The choice of bandwidth will be discussed in detail in Lemma A.11
of the Supplementary Material.

The output of Algorithm \ref{alg1} then combines the debiased estimators of two linear functionals, i.e., $\widehat{\mathbf{x}_{\text {new}}\trans \boldsymbol{\beta}_{\tau,1}}$ and $\widehat{\mathbf{x}_{\text {new}}\trans \boldsymbol{\beta}_{\tau,2}}$, and estimates the IQTE
$\Delta_{\tau,\text {new}}$ by 
$$\widehat{\Delta}_{\tau,\text {new}}=\widehat{\mathbf{x}_{\text {new}}\trans \boldsymbol{\beta}_{\tau,1}}-\widehat{\mathbf{x}_{\text {new}}\trans \boldsymbol{\beta}_{\tau,2}}.\vspace{-2mm}$$ 

\subsection{Statistical Inference for \texorpdfstring{$\Delta_{\tau,\text {new}}$}{DeltaNew}}

Finally, we turn to the inference procedures for $\Delta_{\tau,\text {new}}$, namely, the confidence interval of $\Delta_{\tau,\text {new}}$ and the hypothesis testing of $H_0:\Delta_{\tau,\text {new}}\leq0$ v.s. $H_1:\Delta_{\tau,\text {new}}>0$, for any $\tau\in\mathcal{U}$.

Based on the decomposition of the debiased estimators in Equation \eqref{dec}, we can estimate the variance of $\widehat{\Delta}_{\tau,\text {new}}$ by 
\begin{equation}\label{var_est}
\widehat{V}_{\tau}=\frac{1}{n_1}\tau(1-\tau)\widehat{\mathbf{M}}\trans_{\tau,1}\widehat{\boldsymbol{\Sigma}}_{\tau,1}\widehat{\mathbf{M}}_{\tau,1}+\frac{1}{n_2}\tau(1-\tau)\widehat{\mathbf{M}}\trans_{\tau,2}\widehat{\boldsymbol{\Sigma}}_{\tau,2}\widehat{\mathbf{M}}_{\tau,2}.
\end{equation}
Therefore, for a pre-specified level $\alpha\in (0,1)$, a confidence interval for $\Delta_{\tau,\text {new}}$ with $(1-\alpha)\times100\%$ coverage can be constructed as follows:
\begin{align}
\label{CI}
\operatorname{CI}_{\alpha}=\left(\widehat{\Delta}_{\tau,\text {new}}-z_{\alpha/2}\sqrt{\widehat{V}_{\tau}},\quad \widehat{\Delta}_{\tau,\text {new}}+z_{\alpha/2}\sqrt{\widehat{V}_{\tau}}\right),
\end{align}
where $z_{\alpha/2}$ is the upper {$\alpha/2$} quantile for the standard normal distribution. 
We verify in Section \ref{theory.sec} that $\operatorname{CI}_{\alpha}$ achieves the $(1-\alpha)\times100\%$ coverage asymptotically. The minimax result for the confidence intervals is also established to show the advantage of CI$_{\alpha}$.

In addition, an $\alpha$-level test for $H_0:\Delta_{\tau,\text {new}}\leq0$ can be developed in the following:
\begin{align}
    \phi_{\alpha}=\mathbb{I}\left\{\widehat{\Delta}_{\tau,\text {new}}-z_{\alpha}\sqrt{\widehat{V}_{\tau}}>0\right\}, \label{test}
\end{align}
and we reject the null hypothesis if $\phi_{\alpha}=1$.
The decision rule implies that if we reject the null hypothesis, the first treatment is considered preferable, given that a larger response value indicates a better outcome for the individual.
Corollary \ref{cor1} in Section \ref{asymp.subsec} shows that $\phi_{\alpha}$ controls type I error under level $\alpha$ asymptotically. 
Moreover, the power of the proposed test is obtained, which is essential for analyzing the minimax optimality of the testing procedures.

\section{Theoretical Results}\label{theory.sec}
This section explores theoretical properties of the proposed estimation and inference procedures. In Section \ref{asymp.subsec}, we perform the decomposition and establish the weak convergence of the debiased estimator  $\widehat{\mathbf{x}_{\text {new}}\trans \boldsymbol{\beta}_{\tau,k}}$. This serves as the foundation for the validity of the confidence interval CI$_\alpha$ and the decision rule $\phi_{\alpha}$, as discussed in Section \ref{minimax.subsec}. Additionally, we provide minimax optimality results for the expected length of confidence intervals and the detection boundary for hypothesis testing procedures.

\subsection{Uniform Convergence of the Debiased Estimator}\label{asymp.subsec}
We first collect some notations and sufficient conditions. 
Define the covariance matrix as $\boldsymbol{\Sigma}_k=\mathbb{E}\left(\mathbf{X}_{k, i}\mathbf{X}_{k, i}\trans\right)$ and the precision matrix as $\boldsymbol{\Omega}_k=\boldsymbol{\Sigma}_k^{-1}$ for $k=1,2$.
Denote by $f$ a general density function with support on $\mathbb{R}$, and let $F$ be the corresponding CDF.
Denote by $\overline{\mathcal{U}}$ a compact set satisfying $\mathcal{U} \subsetneq \overline{\mathcal{U}} \subset (0,1)$.

\begin{condition}\label{cond1}
    Suppose that $\mathbf{X}_{k, i} \in \mathbb{R}^p$ is sub-Gaussian, i.e. $\|\mathbf{X}_{k, i}\|_{\psi_2}\leq C_0$ for some constant $C_0>0$, $1\leq i \leq n_k$, $k=1,2$. 
    Suppose that $\boldsymbol{\Sigma}_k\in\Gamma=\{\boldsymbol{\Sigma}\in\mathbb{R}^{p\times p}: c_1 \leq \lambda_{\text{min}}(\boldsymbol{\Sigma})  \leq \lambda_{\text{max}}(\boldsymbol{\Sigma}) \leq C_1\}$ for some constants $C_1 \geq c_1 >0$. 
\end{condition}

\begin{condition}\label{cond2}
    Suppose that  $f_{{Y}_{k,i}\mid\mathbf{X}_{k,i}}(\cdot \vert \mathbf{x}) \in \mathcal{F}=\{f:\vert f^{(l)}(\cdot) \vert \leq C_2,l=0,1, 2, \ f(F^{-1}(\tau))\geq c_2\text{ uniformly in }\tau\in{\overline{\mathcal{U}}}\}$, for any $\mathbf{x}\in \mathbb{R}^p$, $1\leq i \leq n_k$, $k=1,2$, and some constants $C_2 \geq c_2 >0$.
\end{condition}

\begin{condition}\label{cond3}
Assume that $\sup\limits_{\tau \in { \overline{\mathcal{U}}}}\Vert\widehat{\boldsymbol{\beta}}_{\tau,k} \Vert_0^{}\lesssim s_k := \sup\limits_{\tau \in \overline{\mathcal{U}}}\Vert{\boldsymbol{\beta}}_{\tau,k} \Vert_0^{}$, $\sup\limits_{\tau \in \overline{\mathcal{U}}}\Vert\widehat{\boldsymbol{\beta}}_{\tau,k}-\boldsymbol{\beta}_{\tau,k}\Vert_{2}^{}\lesssim \left(\frac{s_k\log p}{n_k}\right)^{1/2}$, with probability larger than $1-g_k(n_k)$, where $g_k(n_k)\rightarrow 0$ as $n_k \rightarrow \infty$, for $1\leq i \leq n_k$ and $k=1,2$.
\end{condition}

\begin{condition}\label{cond4}
    Assume $\lambda_k \asymp \sqrt{\log p/n_k}$ and $\gamma_k \asymp \sqrt{\log p}$, $k=1,2$.
\end{condition}

\begin{remark} \label{remark1}
Condition \ref{cond1} is mild: the sub-Gaussian assumption is commonly assumed in the literatures on treatment effect inference \citep[e.g.,][]{athey2018,wang2020,wu2023,bradic2023highdimensional}; 
the eigenvalue assumption ensures a well-conditioned covariance matrix of the covariates and is widely adopted as well \citep[e.g.,][]{Belloni2011,vandegeer2014,xia2018,he2023}.
Condition \ref{cond2} assumes differentiable and bounded densities for the conditional distribution of responses given covariates, and is frequently employed in high dimensional quantile regression literatures \citep[e.g.,][]{fan2014,bradic2017uniform,giessing2023debiased}. 
Many common distributions meet these requirements, including normal distribution, Student's $t$ distribution, asymmetric Laplace distribution, and many others. 
Condition \ref{cond3} specifies the required convergence rate for the initial coefficient estimator, which is readily achieved using methods such as those proposed by \cite{Belloni2011} and \cite{zhao2015general}, among others. 
Finally, Condition \ref{cond4} presents the orders of the tuning parameters $\lambda_k$ and $\gamma_k$ in the construction of the projection direction.
\end{remark}
\vspace{-4mm}
Recall that we decompose the linear functional estimators $\widehat{\mathbf{x}_{\text {new}}\trans \boldsymbol{\beta}_{\tau,k}}$ as follows:
\begin{align*}
        \widehat{\mathbf{x}_{\text {new}}\trans \boldsymbol{\beta}_{\tau,k}}-\mathbf{x}_{\text {new}}\trans \boldsymbol{\beta}_{\tau,k}=U_{\tau,k}+\Delta_{\tau,k},
    \end{align*}
where the bias term $\Delta_{\tau,k}$ can be further decomposed into five refined components as shown in \eqref{bias.decomp}.
Then the decomposition can be expressed more explicitly by
\begin{align*}
    \widehat{\mathbf{x}_{\text {new}}\trans \boldsymbol{\beta}_{\tau,k}}-\mathbf{x}_{\text {new}}\trans \boldsymbol{\beta}_{\tau,k}=U_{\tau,k}+\Delta_{\tau,k,1}+\Delta_{\tau,k,2}+\Delta_{\tau,k,3}+\Delta_{\tau,k,4}+\Delta_{\tau,k,5}.
\end{align*}

To establish the weak convergence for the proposed debiased estimator $\widehat{\Delta}_{\tau,\text {new}}$, we first demonstrate the performance of the bias terms.
Recall that Constraint \eqref{1} in Algorithm \ref{alg1} is developed to bound $\Delta_{\tau,k,1}$, while Constraint \eqref{3} is used to control the remaining four bias terms. 
Then, based on the convergence rates of the sparsity function estimators established in Lemma A.11
of the Supplementary Material, we derive the convergence rates for all five bias terms, which subsequently lead to the uniform convergence rate for $\Delta_{\tau,k}$, as  established in the following proposition.

\begin{proposition}    \label{prop1}
    Under Conditions \ref{cond1} - \ref{cond4} and assume that $s_k\log p/n_k^{{1}/{3}}=o(1)$, then with probability larger than $1-p^{-\kappa}-g_k(n_k)$ for some constant $\kappa>0$, we have for $k=1,2$,
    \begin{align*}
        \sup _{\tau \in \mathcal{U}}\left|\Delta_{\tau,k}\right|\lesssim\left\|\mathbf{x}_{\text {new}}\right\|_2 \frac{s_k^{7/6}(\log p)^{5/3}}{{n_k}^{5/6}}.
    \end{align*}
\end{proposition}
\begin{remark}
It is worth noting that deriving the convergence rates for the bias terms in the context of quantile regression presents additional challenges compared to the ATE results in \cite{cai2023statistical} under the linear regression setting. 
First, the non-differentiable score function $\varphi(\cdot)$ utilized in the debiasing step \eqref{4} introduces technical difficulties.
Second, unlike focusing solely on the conditional mean, achieving uniform convergence across a set of quantile levels requires more nuanced techniques, where sophisticated empirical process theories and concentration inequalities are employed to attain the goal.
\end{remark}

By Proposition \ref{prop1}, if the sparsity level of $\boldsymbol{\beta}_{\tau,k}$ is not growing too fast, the bias term can be asymptotically negligible. We briefly remark that, beyond the high-dimensional QTE inference addressed in this article, the result in Proposition \ref{prop1} is also of independent interest, especially in the statistical inference on linear functionals, i.e., $\mathbf{x}_{\text {new}}\trans \boldsymbol{\beta}_{\tau,k}$; see for example \cite{cai2017}.

To ensure valid inference, we further provide bounds for the variance estimator of $\widehat{\Delta}_{\tau,\text {new}}$ as defined in \eqref{var_est}.
\begin{proposition}\label{prop2}
    Under Conditions \ref{cond1} - \ref{cond4} and assume that $ s_k\log p/\sqrt{n_k}=o(1)$, then with probability larger than $1-p^{-\kappa}-g_k(n_k)$ for some constant $\kappa>0$, it holds that 
    $$\sup\limits_{\tau\in\mathcal{U}}\sqrt{\widehat{V}_\tau}\asymp \left\|\mathbf{x}_{\text {new}}\right\|_2\left(\frac{1}{\sqrt{n_1}}+\frac{1}{\sqrt{n_2}}\right).$$
\end{proposition}
\begin{remark}\label{remark.structure-free}
Note that Constraint \eqref{2} in Algorithm \ref{alg1} is crucial for establishing the lower bound of the variance estimator with structure-free covariates $\mathbf{x}_{\text {new}}$. Intuitively, without Constraint \eqref{2}, the vector $\widehat{\mathbf{M}}_{\tau,k}=\mathbf{0}$ is included in the feasible set of the optimization problem \eqref{obj}, provided that $\Vert \mathbf{x}_{\text {new}} \Vert_{\infty}^{} \leq \Vert \mathbf{x}_{\text {new}} \Vert_{2}^{}\lambda_k$. Consequently, the variance can no longer dominate the bias for such dense loading vectors, leading to a failure of the debiasing procedure and the subsequent inference methods.
\end{remark}

Combining the above two propositions, we next present the main result: the debiased estimator $\widehat{\mathbf{x}_{\text {new}}\trans \boldsymbol{\beta}_{\tau,k}}$ can be approximated by a Gaussian process on $\mathcal{U}$. To facilitate the presentation of this result, we first construct an auxiliary projection vector $\widetilde{\mathbf{M}}_{\tau,k}$ by the following procedure:
\begin{align}
    \widetilde{\mathbf{M}}_{\tau,k}=\underset{\mathbf{M}_{\tau,k} \in \mathbb{R}^p}{\arg \min\ } \mathbf{M}_{\tau,k}\trans \widetilde{\boldsymbol{\Sigma}}_{\tau,k} \mathbf{M}_{\tau,k}, \text { s.t. } &\Vert\widehat{\boldsymbol{\Sigma}}_k \mathbf{M}_{\tau,k}-\mathbf{x}_{\text {new}}\Vert_{\infty} \leq\left\|\mathbf{x}_{\text {new}}\right\|_2 \lambda_k, \nonumber\\
    & |\mathbf{x}_{\text {new}}\trans \widehat{\boldsymbol{\Sigma}}_k \mathbf{M}_{\tau,k}-\left\|\mathbf{x}_{\text {new}}\right\|_2^2| \leq\left\|\mathbf{x}_{\text {new}}\right\|_2^2 \lambda_k,\nonumber\\
    &\underset{1\leq i\leq n_k}{\max} |\mathbf{X}_{k, i}\trans\mathbf{M}_{\tau,k} | \leq \left\|\mathbf{x}_{\text {new}}\right\|_2\gamma_k, \label{auxi}
\end{align}
where $\widetilde{\boldsymbol{\Sigma}}_{\tau,k}=n_k^{-1}\sum\limits_{i=1}^{n_k} \left({1}/{f_{k,i}(F_{k,i}^{-1}(\tau))} \right)^2\mathbf{X}_{k, i} \mathbf{X}_{k, i}\trans$. 
It can be seen that the only difference between $\widetilde{\mathbf{M}}_{\tau,k}$ and $\widehat{\mathbf{M}}_{\tau,k}$ lies in the objective function of the constrained optimization problem \eqref{auxi}. By such construction, the randomness in $\widetilde{\mathbf{M}}_{\tau,k}$ solely originates from $\mathbf{X}_k$. 

Based on this auxiliary projection,  for any $\tau_1,\tau_2\in\mathcal{U}$, we define
$$\widetilde{V}_{\tau,k}={\widetilde{\mathbf{M}}_{\tau,k}\trans\widetilde{\Sigma}_{\tau,k}\widetilde{\mathbf{M}}_{\tau,k}},\ \widetilde{J}_k(\tau_1,\tau_2)=\frac{1}{n_k}\sum\limits_{i=1}^{n_k}\frac{1}{f_{k,i}(F^{-1}_{k,i}(\tau_1))f_{k,i}(F^{-1}_{k,i}(\tau_2))}\widetilde{\mathbf{M}}_{\tau_1,k}\trans\mathbf{X}_{k,i}\mathbf{X}_{k,i}\trans\widetilde{\mathbf{M}}_{\tau_2,k},$$
which respectively approximate the variance of $U_{\tau,k}$ and the covariance between $U_{\tau_1,k}$ and $U_{\tau_2,k}$ up to a constant factor.
Then, the weak convergence result is established in the following theorem.
\begin{theorem}[Weak convergence]\label{theowc}
    Under Conditions \ref{cond1} - \ref{cond4} and assume that 
    $s_k=o(n_k^{\frac{2}{7}}/\left(\log p\right)^{\frac{10}{7}})$, then
    $${\sqrt{n_1}\widetilde{V}_{\tau,1}^{-\frac{1}{2}}}\left(\widehat{\mathbf{x}_{\text {new}}\trans \boldsymbol{\beta}_{\tau,1}}-\mathbf{x}_{\text {new}}\trans \boldsymbol{\beta}_{\tau,1}\right)+{\sqrt{n_2}\widetilde{V}_{\tau,2}^{-\frac{1}{2}}}\left(\widehat{\mathbf{x}_{\text {new}}\trans \boldsymbol{\beta}_{\tau,2}}-\mathbf{x}_{\text {new}}\trans \boldsymbol{\beta}_{\tau,2}\right)=\mathbb{G}_{1}(\tau)+\mathbb{G}_{2}(\tau)+e(\tau),$$
    where $\mathbb{G}_1(\cdot)$ and $\mathbb{G}_2(\cdot)$ are two independent centered Gaussian processes on $\mathcal{U}$, conditional on $(\mathbf{X}_1,\mathbf{X}_2)$, with covariance functions between $\tau_1,\tau_2 \in \mathcal{U}$ given by
    \begin{align}
        \widetilde{W}_{n_k}(\tau_1,\tau_2)=\left(\tau_1 \wedge \tau_2-\tau_1 \tau_2\right)\widetilde{V}_{\tau_1,k}^{{-\frac{1}{2}}}\widetilde{J}_k(\tau_1,\tau_2)\widetilde{V}_{\tau_2,k}^{{-\frac{1}{2}}}, \nonumber
    \end{align}
    and the remaining error satisfies that $\sup\limits_{\tau\in\mathcal{U}}|e(\tau)|=o_P(1)$, as $n_1\wedge n_2 =:n \rightarrow \infty$. 
\end{theorem}

\begin{remark}
We briefly discuss the conditions and technical challenges involved.
First, we impose no structural assumptions on the loading vector $\mathbf{x}_{\text {new}}$, unlike the $\ell_1$-convergent or $\ell_2$-bounded conditions in \cite{bradic2017uniform,yan2023,giessing2023debiased}.
Second, as the dimension $p$ increases with the sample size $n_k$, we establish uniform convergence over these so-called changing classes \citep{vaart1996}, where classical Donsker theorems for the weak convergence of finite dimensional quantile processes  \citep{koenker2002process,angrist2006} are not applicable.
Third, the technical tool utilized here is based on Gaussian coupling, which is fundamentally different from the approaches used in existing works on the quantile process of linear functionals \citep{chao2017,chao2019,giessing2023debiased}.
\end{remark}

\subsection{Inference Validity and Minimax Optimality}\label{minimax.subsec}
Next, we turn to the inference procedures of $\Delta_{\tau,\text {new}}$ for any fixed $\tau\in \mathcal{U}$. 
Combining Proposition \ref{prop2} and Theorem \ref{theowc}, we attain the asymptotic normality result for the proposed debiased estimator of ${\Delta}_{\tau,\text {new}}$ in the following corollary. 
\begin{corollary}[Asymptotic normality]\label{theo1}
    Under the same conditions of Theorem \ref{theowc},
    for any fixed $\tau \in \mathcal{U}$, we have
    $$\frac{1}{{\sqrt{\widehat{V}_{\tau}}}}\left(\widehat{\Delta}_{\tau,\text {new}}-{\Delta}_{\tau,\text {new}}\right)\xrightarrow[]{d}N(0,1), \mbox{ as } n \rightarrow \infty.$$
\end{corollary}

Based on Corollary \ref{theo1}, we further obtain the following two results on the proposed confidence interval CI$_{\alpha}$, as defined in \eqref{CI}, and the proposed decision rule $\phi_{\alpha}$, as defined in \eqref{test}, for testing $H_0:\Delta_{\tau,\text {new}}\leq0$.

\begin{corollary}\label{cor1}
    Under the same conditions in Corollary \ref{theo1}, we have
    $$
    \liminf_{n \rightarrow \infty} \mathbb{P}\left(\Delta_{\tau,\text {new}} \in \mathrm{CI}_\alpha\right) \geq 1-\alpha ,\text{ for any } \mathbf{x}_{\text {new}} \in \mathbb{R}^p.
    $$
\end{corollary}

\begin{corollary}\label{cor2}
    Under the same conditions in Corollary \ref{theo1}, we have
    $$\limsup_{n \rightarrow \infty} \mathbb{P}_{H_0}\left(\phi_{\alpha}=1\right) \leq \alpha ,\text{ for any } \mathbf{x}_{\text {new}} \in \mathbb{R}^p.$$
\end{corollary}

The two corollaries demonstrate the asymptotic $(1-\alpha)\times100\%$ coverage for the proposed confidence interval, as well as the asymptotic type I error control for the proposed testing procedure, respectively.

\begin{remark}\label{x_new}
It is important to remark that, the results in Corollaries \ref{theo1} - \ref{cor2} do not rely on any structural assumptions on the individual loading vector $\mathbf{x}_{\text {new}}$, thanks to the introduction of  Constraint \eqref{2}. Without such constraint, ensuring the validity of the inference procedures would require more stringent assumptions on the loading vector, such as $\ell_2$-convergence or boundedness \citep[e.g.,][]{zhao2015general, bradic2017uniform, giessing2023debiased}.
\end{remark}

To demonstrate the superiority of our proposed debiased estimator, we next establish the minimax optimality results for both confidence intervals and hypothesis testing procedures. 
Denote by $f_{k}(\cdot\vert \mathbf{x}):=f_{Y_{k,i}\vert \mathbf{X}_{k,i}}(\cdot\vert \mathbf{x})$ for $k=1,2$ and consider the parameter space $\Theta$ defined by 
$$
    \Theta= \left\{ \boldsymbol{\theta} = \{\boldsymbol{\beta}_{\tau,k}, \boldsymbol{{\Sigma}}_k,f_{k}(\cdot\vert \mathbf{x})\}_{k=1,2}:\|\boldsymbol{\beta}_{\tau,k}\|_0^{} \leq s_k,\boldsymbol{{\Sigma}}_k\in\Gamma,f_{k}(\cdot\vert \mathbf{x})\in\mathcal{F} \text{ for any }\mathbf{x}\in\mathbb{R}^p\right\}.
$$

We start with the minimax optimal framework for confidence intervals.
For $0<\alpha<1$, let $\mathcal{I}_\alpha(\Theta;\mathbf{x}_{\text {new}}, \mathbf{Z})$ denote the set of all confidence intervals for $\Delta_{\tau,\text {new}}$ that have an asymptotic $(1-\alpha)$ coverage probability over the parameter space $\Theta$,  constructed based on the observed data $\mathbf{Z}=\left(\mathbf{Y}_1, \mathbf{X}_1, \mathbf{Y}_2, \mathbf{X}_2\right)$, i.e., for any $\mathbf{x}_{\text {new}}\in \mathbb{R}^p$ we define 
\begin{align}
    \mathcal{I}_\alpha(\Theta;\mathbf{x}_{\text {new}}, \mathbf{Z})&=\left\{\left(l(\mathbf{x}_{\text {new}}, \mathbf{Z}), u(\mathbf{x}_{\text {new}}, \mathbf{Z})\right):\vphantom{}\right.\nonumber\\
&\quad\quad\quad\left.\inf _{\boldsymbol{\theta} \in \Theta} \mathbb{P}_{\boldsymbol{\theta}}(l(\mathbf{x}_{\text {new}}, \mathbf{Z}) \leq \Delta_{\tau,\text {new}} \leq u(\mathbf{x}_{\text {new}}, \mathbf{Z})) \geq 1-\alpha+o(1), \mbox{ as } n \rightarrow \infty\right\},\nonumber
\end{align}
where $l(\mathbf{x}_{\text {new}}, \mathbf{Z}), u(\mathbf{x}_{\text {new}}, \mathbf{Z}) \in \mathbb{R}$ are the lower and upper bounds of the confidence intervals constructed from $\mathbf{Z}$.
For any valid confidence interval $\mathrm{CI}_\alpha(\mathbf{x}_{\text {new}}, \mathbf{Z})\in \mathcal{I}_\alpha(\Theta;\mathbf{x}_{\text {new}}, \mathbf{Z})$, we evaluate its performance by measuring its maximal expected length over $\Theta$, i.e.
$$
L\left\{\Theta;\mathrm{CI}_\alpha(\mathbf{x}_{\text {new}}, \mathbf{Z})\right\}=\sup _{\boldsymbol{\theta} \in \Theta} \mathbb{E}_{\boldsymbol{\theta}} L\left\{ \mathrm{CI}_\alpha(\mathbf{x}_{\text {new}}, \mathbf{Z})\right\},
$$
where $L\left\{ \mathrm{CI}_\alpha(\mathbf{x}_{\text {new}}, \mathbf{Z})\right\}=u(\mathbf{x}_{\text {new}}, \mathbf{Z})-l(\mathbf{x}_{\text {new}}, \mathbf{Z})$ is the length of the interval. 
In addition, we define $L_\alpha^*\left(\Theta;\mathbf{x}_{\text {new}}\right)$ as the infimum of the maximal expected length among all asymptotic $(1-\alpha)$ level confidence intervals over $\Theta$, namely,
$$
L_\alpha^*\left(\Theta;\mathbf{x}_{\text {new}}\right)=\inf _{\mathrm{CI}_\alpha(\mathbf{x}_{\text {new}}, \mathbf{Z}) \in \mathcal{I}_\alpha(\Theta;\mathbf{x}_{\text {new}}, \mathbf{Z})} L\left\{ \Theta;\mathrm{CI}_\alpha(\mathbf{x}_{\text {new}}, \mathbf{Z})\right\}.
$$
This quantity serves as a benchmark to measure whether a confidence interval is minimax rate-optimal. Specifically, we say $\mathrm{CI}_\alpha(\mathbf{x}_{\text {new}}, \mathbf{Z})\in \mathcal{I}_\alpha(\Theta;\mathbf{x}_{\text {new}}, \mathbf{Z})$ attains minimax rate optimality asymptotically, if it satisfies
$$\lim_{n \rightarrow \infty} \mathbb{P}\left(L\left\{ \mathrm{CI}_\alpha(\mathbf{x}_{\text {new}}, \mathbf{Z})\right\} \asymp L_\alpha^*\left(\Theta;\mathbf{x}_{\text {new}}\right)\right)=1.$$ 

The following theorem presents the minimax expected length for the confidence intervals with guaranteed coverage over the parameter space defined above, and  shows that our proposed CI$_\alpha$ defined in \eqref{CI}
meets the optimality asymptotically. 

\begin{theorem}[Minimaxity of CI$_\alpha$]\label{theo2}
    Under the conditions of Corollary \ref{theo1}, if $0<\alpha<\frac{1}{2}$, $\max\left\{s_1, s_2\right\}\leq \frac{p}{4}$, then for any $ \mathbf{x}_{\text {new}}\in \mathbb{R}^p$ with $\Vert\mathbf{x}_{\text {new}}\Vert_0^{}=q\leq 2\left(s_1+s_2\right)$, it holds  that
    $$L_\alpha^*\left(\Theta;\mathbf{x}_{\text {new}}\right) \asymp \|\mathbf{x}_{\text {new}}\|_2^{}\left(\frac{1}{\sqrt{n_1}}+\frac{1}{\sqrt{n_2}}\right).$$
    Furthermore, $\mathrm{CI}_\alpha$ is asymptotically minimax rate-optimal, namely,  
    $$\mathrm{CI}_\alpha \in \mathcal{I}_\alpha(\Theta;\mathbf{x}_{\text {new}}, \mathbf{Z})\text{ and }\lim_{n \rightarrow \infty} \mathbb{P}\left(L\left\{\mathrm{CI}_\alpha\right\}\asymp L_\alpha^*\left(\Theta;\mathbf{x}_{\text {new}}\right)\right)=1.$$
\end{theorem}
\begin{remark}
The additional conditions on $s_k$ and $\Vert\mathbf{x}_{\text {new}}\Vert_0^{}$ are crucial for the least favorable construction in obtaining the lower bound of $L_\alpha^*\left(\Theta;\mathbf{x}_{\text {new}}\right)$. They are more relaxed compared to the linear regression setting in \cite{cai2017}, which assumes $s_k\leq p^c$ for some $0<c<1/2$.
\end{remark}

Parallel to the minimaxity of confidence intervals, we next develop an optimality framework for hypothesis testing of the IQTE ${\Delta}_{\tau,\text {new}}$. 
For a given covariate $\mathbf{x}_{\text {new}}\in\mathbb{R}^p$, let $\mathcal{H}_0(\Theta;\mathbf{x}_{\text {new}})=\left\{{\boldsymbol{\theta}}\in\Theta:{\Delta}_{\tau,\text {new}}\leq0\right\}$ denote the null hypothesis space and $\mathcal{H}_1(\delta,\Theta;\mathbf{x}_{\text {new}})=\left\{{\boldsymbol{\theta}}\in\Theta:{\Delta}_{\tau,\text {new}}=\delta\right\}$ represent a simple alternative hypothesis space with detection boundary {$\delta>0$}.
Then for a test function $\phi(\mathbf{x}_{\text {new}}, \mathbf{Z})$, its type I error and power are respectively defined by 
\begin{align*}
\boldsymbol{\alpha}(\phi(\mathbf{x}_{\text {new}}, \mathbf{Z}); \mathcal{H}_0(\Theta;\mathbf{x}_{\text {new}}))&=\sup _{\theta \in \mathcal{H}_0(\Theta;\mathbf{x}_{\text {new}})} \mathbb{E}_\theta \phi(\mathbf{x}_{\text {new}}, \mathbf{Z}),\\
\boldsymbol{\omega}(\phi(\mathbf{x}_{\text {new}}, \mathbf{Z}); \mathcal{H}_1(\delta,\Theta;\mathbf{x}_{\text {new}}))&=\inf _{\theta \in \mathcal{H}_1(\delta,\Theta;\mathbf{x}_{\text {new}})} \mathbb{E}_\theta \phi(\mathbf{x}_{\text {new}}, \mathbf{Z}).
\end{align*}
Intuitively, it is easier to differentiate null and alternative spaces if the detection boundary $\delta$ is larger. 
Therefore, 
we aim to find the smallest detection boundary $\delta$ among all asymptotic $\alpha$-level test with maximal power exceeding $1-w$ asymptotically for some constant $w \in \left(0, 1\right)$, i.e., 
$$
\delta^{*}\left(\Theta;\mathbf{x}_{\text {new}}\right)=\underset{\delta}{\arg \inf }\left\{\delta: \sup _{\phi: \boldsymbol{\alpha}(\phi) \leq \alpha+o(1)} \boldsymbol{\omega}(\phi; \mathcal{H}_1(\delta,\Theta;\mathbf{x}_{\text {new}})) \geq 1-w+o(1), \mbox{ as }n\rightarrow \infty\right\},
$$
where we abbreviate $\phi(\mathbf{x}_{\text {new}}, \mathbf{Z})$ and $\boldsymbol{\alpha}(\phi(\mathbf{x}_{\text {new}}, \mathbf{Z}); \mathcal{H}_0(\Theta;\mathbf{x}_{\text {new}}))$ to $\phi$ and $\boldsymbol{\alpha}(\phi)$, respectively, when there is no ambiguity. The above quantity serves as a benchmark to measure whether a test function is minimax rate-optimal. Specifically, if a test $\phi(\mathbf{x}_{\text {new}}, \mathbf{Z})$ satisfies 
\begin{align}
    \boldsymbol{\alpha}(\phi(\mathbf{x}_{\text {new}}, \mathbf{Z}); \mathcal{H}_0(\Theta;\mathbf{x}_{\text {new}})) \leq \alpha+o(1),  \label{mini1}
\end{align}
and 
\begin{align}
    \boldsymbol{\omega}(\phi(\mathbf{x}_{\text {new}}, \mathbf{Z}); \mathcal{H}_1(\delta,\Theta;\mathbf{x}_{\text {new}})) \geq 1-w+o(1) \text { for } \delta \gtrsim \delta^{*}\left(\Theta;\mathbf{x}_{\text {new}}\right),\label{mini}
\end{align}
then $\phi(\mathbf{x}_{\text {new}}, \mathbf{Z})$ is an asymptotically optimal test function.

The following theorem establishes the minimax rate for the detection boundary of the hypothesis testing procedures with asymptotic type I error control.
In addition, we provide a power analysis of the proposed test $\phi_\alpha$ in \eqref{test}, which demonstrates that $\phi_\alpha$ achieves asymptotic minimaxity.

\begin{theorem}[Minimaxity of $\phi_\alpha$]\label{theo3}
    Under the conditions of Corollary \ref{theo1}, if $w<1-\alpha$, then for any $\mathbf{x}_{\text {new}}$ with $\Vert\mathbf{x}_{\text {new}}\Vert_0^{}=q\leq (s_1\wedge s_2)$,
    $$
    \delta^{*}\left(\Theta;\mathbf{x}_{\text {new}}\right) \asymp \left\|\mathbf{x}_{\text {new}}\right\|_{2}\left(\frac{1}{\sqrt{n_1}}+\frac{1}{\sqrt{n_2}}\right).
    $$
    In addition, let $\delta=K/(\sqrt{c_1}c_2)\sqrt{\tau(1-\tau)}\left\|\mathbf{x}_{\text {new}}\right\|_2\left(\frac{1}{\sqrt{n_1}}+\frac{1}{\sqrt{n_2}}\right)$ where $K>0$ is a constant, then for any $\mathbf{x}_{\text {new}} \in \mathbb{R}^p$, 
   $$\liminf_{n \rightarrow \infty} \boldsymbol{\omega}(\phi_{\alpha}; \mathcal{H}_1(\delta,\Theta;\mathbf{x}_{\text {new}}))\geq 1-\Phi(z_{\alpha}-K),$$
   where $\Phi(\cdot)$ is the CDF of the standard normal distribution.
    Consequently, by taking $K\geq z_{\alpha}-\Phi^{-1}\left(w\right)$, $\phi_\alpha$ satisfies \eqref{mini1} and \eqref{mini} and is minimax rate-optimal.
\end{theorem}
\begin{remark}
The minimax lower bounds for  two-sample IQTE, established in Theorems \ref{theo2} and \ref{theo3}, require more advanced techniques than those used for one-sample IATE confidence intervals \citep{cai2017} and hypothesis testing \citep{cai2021}.   
   Specifically, we relax the Gaussian assumptions on $Y_{k,i}\vert \mathbf{X}_{k,i}$ and accommodate heavy-tailed distributions as long as Condition \ref{cond2} is satisfied. Therefore, the technical tools used under normality \citep{Cai_2012} are not directly applicable, making the derivation of the minimax lower bound particularly challenging. We address this by leveraging theoretical results related to the asymmetric Laplace distribution to achieve the desired minimax optimality. As a consequence, our approach allows $\|\mathbf{x}_{\text {new}}\|_0^{}$ to be of the same order as the dimension of the covariates $p$ for the lower bound derivation, accommodating much broader settings compared to the ATE framework. 
\end{remark}

\section{Simulation}\label{simu.sec}
In this section, we demonstrate the numerical performance of our proposed estimator $\widehat{\Delta}_{\tau,\text {new}}$ (denoted by \texttt{IQTE}) and compare it with a few competing methods: the  lasso-type method \citep{Belloni2011} which directly employs the regularized estimators in \eqref{Belloni.est} without debiasing (denoted by \texttt{Lasso}), the plug-in debiased estimator \citep[e.g.,][]{zhao2015general} as calculated by Equation \eqref{naive.est} (denoted by \texttt{Deb}), and the regression rank-score debiased estimator \citep{giessing2023debiased} (denoted by \texttt{RS}).
We first compare the biases and standard errors of the aforementioned estimators. Next, the coverage probabilities and lengths of the confidence intervals are calculated to investigate the validity and efficiency of the corresponding procedures. Finally, the type I error and empirical power of the testing procedures based on the above estimators are collected and compared.

\subsection{Data Generation and Implementation}
We set the dimensionality of the covariates $p=600$, and consider a range of choices for the sample sizes of two treatment groups, i.e., $(n_1,n_2)\in\left\{(300,300),(500,550),(900,800)\right\}$. 
For $1\leq i \leq n_k$, $k=1,2$, we independently generate $\widetilde{\mathbf{X}}_{k,i}$ from a multivariate Gaussian distribution with mean zero and covariance matrix $\boldsymbol{\Sigma}_k=\left({\sigma}_{ij}\right)_{i,j=1}^{p}$, where $\sigma_{ij}=0.5^{|i-j|}$ for $1\leq i,j\leq p$.
Next, the responses for two treatment groups are generated by
\begin{align}
  {Y}_{1,i}&=\sum\limits_{j=3}^{p}(\widetilde{\mathbf{X}}_{1,i})_j(\boldsymbol{\beta}_{\tau,1})_{j}+{\left( \vert (\widetilde{\mathbf{X}}_{1,i})_1\vert + 2 \vert (\widetilde{\mathbf{X}}_{1,i})_2 \vert \right){\epsilon}_{1,i}},\text{ for }1\leq i\leq n_1,\label{simu1}\\
  {Y}_{2,i}&=\sum\limits_{j=3}^{p}(\widetilde{\mathbf{X}}_{2,i})_j(\boldsymbol{\beta}_{\tau,2})_{j}+{\left( 2\vert (\widetilde{\mathbf{X}}_{2,i})_1 \vert +  \vert (\widetilde{\mathbf{X}}_{2,i})_2 \vert \right){\epsilon}_{2,i}},\text{ for }1\leq i\leq n_2,\label{simu2}
\end{align}
where $(\boldsymbol{\beta}_{\tau,1})_{3}^{}=2,(\boldsymbol{\beta}_{\tau,1})_{4}^{}=1.6,(\boldsymbol{\beta}_{\tau,1})_{j}=-0.8(j-4)\cdot\mathbb{I}\left\{5\leq j\leq11\right\}$ for $j=5,\dots,p$, $(\boldsymbol{\beta}_{\tau,2})_{3}^{}=1.2,(\boldsymbol{\beta}_{\tau,2})_{4}^{}=0.5,(\boldsymbol{\beta}_{\tau,2})_{j}=0.5(j+1)\cdot\mathbb{I}\left\{5\leq j\leq11\right\}$ for $j=5,\dots,p$, with $(\boldsymbol{\boldsymbol{\xi}})_{j}$ representing the $j{\text {-th}}$ element of any vector $\boldsymbol{\boldsymbol{\xi}}\in\mathbb{R}^p$, and the
random errors $\{\epsilon_{k,i}\}$ are generated independently and identically from the standard normal distribution.  Let 
${\mathbf{X}}_{k}=(\vert \widetilde{\mathbf{X}}_{k,1} \vert,\vert \widetilde{\mathbf{X}}_{k,2} \vert,\widetilde{\mathbf{X}}_{k,3},\ldots,\widetilde{\mathbf{X}}_{k,p})\trans$ for $k=1,2$, then Equations \eqref{simu1} and \eqref{simu2} can be equivalently formulated as
\begin{align}
  F^{-1}_{1,i}(\tau)&=\mathbf{X}_{1,i}\trans\left(z_{1-\tau},2z_{1-\tau},(\boldsymbol{\beta}_{\tau,1})_{3}^{},\ldots,(\boldsymbol{\beta}_{\tau,1})_{p}\right)\trans,\label{simu11}\\
  F^{-1}_{2,i}(\tau)&=\mathbf{X}_{2,i}\trans\left(2z_{1-\tau},z_{1-\tau},(\boldsymbol{\beta}_{\tau,2})_{3}^{},\ldots,(\boldsymbol{\beta}_{\tau,2})_{p}\right)\trans .\label{simu22}
\end{align}
Next, we consider both $\ell_2$-dense and $\ell_2$-sparse loading vectors $\mathbf{x}_{\text{new}}$ as follows:
\begin{setting}[\textbf{Dense loading}]\label{set1}
  Generate $\mathbf{x}_{\text{initial}} \sim N(0,\text{I}_p)$, let $(\mathbf{x}_{\text{new}})_j=0.75\cdot\mathbb{I}\left\{j=1\right\}+5.25\cdot\mathbb{I}\left\{j=2\right\}-1.5\cdot\mathbb{I}\left\{j=3\right\}+1.05\cdot\mathbb{I}\left\{j=4\right\}+\mathbb{I}\left\{j\geq 12\right\}\cdot(\mathbf{x}_{\text{initial}})_j$, for $j=1,\dots,p.$
\end{setting}
\begin{setting}[\textbf{Sparse loading}]\label{set2}
  Generate $\mathbf{x}_{\text{initial}} \sim N(0,\text{I}_p)$, let $(\mathbf{x}_{\text{new}})_j=1.5\cdot\mathbb{I}\left\{j=1\right\}+0.75\cdot\mathbb{I}\left\{j=2\right\}+(-1)^{\mathbb{I}\left\{j=5\right\}}\cdot0.125\cdot(9-j)\cdot\mathbb{I}\left\{3\leq j\leq7\right\}+\mathbb{I}\left\{12\leq j\leq100\right\}\cdot(\mathbf{x}_{\text{initial}})_j$ for $j=1,\dots,p.$
\end{setting}

In order to investigate the quantile heterogeneity, three quantile levels are considered: $\tau_1=0.2,\ \tau_2=0.5\text{, and }\tau_3=0.8$. Under the circumstances above, it can be calculated that ${\Delta}_{\tau_1,\text {new}}=-3.832$, ${\Delta}_{\tau_2,\text {new}}=-0.045$, ${\Delta}_{\tau_3,\text {new}}=3.742$ for Setting \ref{set1}, and ${\Delta}_{\tau_1,\text {new}}=1.906$, ${\Delta}_{\tau_2,\text {new}}=1.275$, ${\Delta}_{\tau_3,\text {new}}=0.644$ for Setting \ref{set2}.  
We set $\alpha = 5\%$ throughout.

To implement our method, the initial coefficient estimator $\widehat{\boldsymbol{\beta}}_{\tau,k}$ is obtained through the R package \texttt{quantreg}, and the bandwidth is set at $h_k=\min\{n_k^{-1/6},\tau(1-\tau)/2\}$, $k=1,2$, to obtain the initial sparsity function estimators $\widehat{\eta}_{k,i}(\tau)$ as suggested by \cite{belloni2019}. Finally, we adopt the five-fold cross-validation to construct the projection direction $\widehat{\mathbf{M}}_{\tau,k}$ in order to reach the proposed debiased estimator $\widehat{\Delta}_{\tau,\text {new}}$.
The comparisons under heavy-tailed distributions are collected in Section B
of the Supplementary Material.
\subsection{Simulation Results} \label{simu2.sec}
We compare the performance of the four competing methods in the following aspects.
First, we present the biases and standard errors of the point estimates for $\Delta_{\tau,\text {new}}$ based on all four approaches (\texttt{IQTE}, \texttt{RS}, \texttt{Deb}, \texttt{Lasso}), with these values evaluated empirically from 1,000 Monte Carlo trials.
Second, we assess the validity of the confidence intervals constructed using \texttt{IQTE}, \texttt{RS} and \texttt{Deb} by comparing their coverage probabilities. 
Note that the \texttt{Lasso} method is excluded from this analysis due to its substantial bias, which undermines its suitability for statistical inference.
Finally, to evaluate the efficiency of the two most competitive methods (\texttt{IQTE} and \texttt{RS}), we compare the lengths of their confidence intervals and the rejection rates of their testing procedures, calculated as the proportion of rejected null hypotheses over 1,000 replications in both settings. 
Note that, under the alternative with $\Delta_{\tau,\text{new}}>0$, the rejection rate reflects the power of the testing procedures, while under the null with $\Delta_{\tau,\text{new}}\leq0$, it indicates the type I error.

\begin{table}[t!] 
\centering
\caption{Comparisons of various methods under dense and sparse settings; $\alpha=5\%$.}
\label{tab1}
\resizebox{\textwidth}{!}{
    \begin{tabular}{cccccccccccccccccc}
      \hline 
      \multicolumn{18}{c}{\textbf{Dense Loading (Setting \ref{set1})}}\\
      \hline 
      & & & \multicolumn{2}{c}{\textbf{Rejection Rate}} & \multicolumn{3}{c}{\textbf{Coverage Rate}} & \multicolumn{2}{c}{\textbf{Length}} & \multicolumn{2}{c}{\textbf{IQTE}} & \multicolumn{2}{c}{\textbf{RS}} & \multicolumn{2}{c}{\textbf{Deb}} & \multicolumn{2}{c}{\textbf{Lasso}} \\  
      \cmidrule(l){4-5} \cmidrule(l){6-8} \cmidrule(l){9-10} \cmidrule(l){11-12} \cmidrule(l){13-14} \cmidrule(l){15-16} \cmidrule(l){17-18} 
      ($n_1$,$n_2$) & $\tau$ & $\Delta_{\tau,\text{new}}$ & IQTE & RS & IQTE & RS & Deb & IQTE & RS & Bias & SE & Bias & SE & Bias & SE & Bias & SE \\
      \hline 
      \multirow[t]{3}{*}{(300,300)} & 0.2 & $\leq0$ & 0.004 & 0.008 & 0.957 & 0.926 & 0.796 & 10.344 & 10.726 & -0.258 & 2.639 & 0.545 & 2.736 & 0.858 & 2.846 & 0.998 & 2.268 \\
      & 0.5 & $\leq0$ & 0.038 & 0.051 & 0.961 & 0.939 & 0.818 & 6.342 & 6.347 & -0.355 & 1.618 & -0.457 & 1.619 & -0.806 & 2.430 & 0.833 & 1.589 \\
      & 0.8 & $>0$ & 0.383 & 0.240 & 0.954 & 0.925 & 0.787 & 10.317 & 10.691 &-0.652 & 2.632 & -0.741 & 2.728 & 0.831 & 2.625 & 0.926 & 2.444 \\
      \hline 
      \multirow[t]{3}{*}{(500,550)} & 0.2 & $\leq0$ & 0.000 & 0.006 & 0.946 & 0.928 & 0.825 & 8.214 & 8.448 & -0.132 & 2.095 & 0.443 & 2.155 & 0.693 & 2.427 & 0.890 & 1.943  \\
      & 0.5 & $\leq0$ & 0.044 & 0.053 & 0.953 & 0.942 & 0.842 & 5.377 & 5.951 & -0.274 & 1.372 & -0.376 & 1.581 & 0.681 & 1.933 & 0.776 & 1.196 \\
      & 0.8 & $>0$ & 0.478 & 0.351 & 0.948 & 0.929 & 0.830 & 8.206 & 8.276 & -0.503 & 2.093 & -0.559 & 2.112 & 0.636 & 2.211 & 0.676 & 1.903 \\   
      \hline 
      \multirow[t]{3}{*}{(900,800)} & 0.2 & $\leq0$ & 0.000 & 0.001 & 0.947 & 0.935 & 0.857 & 6.371 & 6.602 & -0.019 & 1.625 & 0.306 & 1.684 & 0.565 & 1.847 & 0.689 & 1.311  \\
      & 0.5 & $\leq0$ & 0.033 & 0.049 & 0.955 & 0.946 & 0.855 & 4.426 & 5.065 & -0.170 & 1.129 & -0.265 & 1.292 & 0.639 & 1.781 & 0.708 & 1.340 \\
      & 0.8 & $>0$ & 0.669 & 0.523 & 0.943 & 0.939 & 0.863 & 6.345 & 7.072 & -0.215 & 1.619 & -0.504 & 1.804 & 0.592 & 1.910 & 0.641 & 1.310 \\
      \hline
      \multicolumn{18}{c}{\textbf{Sparse Loading (Setting \ref{set2})}}\\
      \hline
      & & & \multicolumn{2}{c}{\textbf{Rejection Rate}} & \multicolumn{3}{c}{\textbf{Coverage Rate}} & \multicolumn{2}{c}{\textbf{Length}} & \multicolumn{2}{c}{\textbf{IQTE}} & \multicolumn{2}{c}{\textbf{RS}} & \multicolumn{2}{c}{\textbf{Deb}} & \multicolumn{2}{c}{\textbf{Lasso}} \\  
      \cmidrule(l){4-5} \cmidrule(l){6-8} \cmidrule(l){9-10} \cmidrule(l){11-12} \cmidrule(l){13-14} \cmidrule(l){15-16} \cmidrule(l){17-18}  
      ($n_1$,$n_2$) & $\tau$ & $\Delta_{\tau,\text{new}}$ & IQTE & RS & IQTE & RS & Deb & IQTE & RS & Bias & SE & Bias & SE & Bias & SE & Bias & SE \\
      \hline
      \multirow[t]{3}{*}{(300,300)} & 0.2 & $>0$ & 0.673 & 0.571 & 0.947 & 0.944 & 0.901 & 3.274 & 4.210 & -0.055 & 0.835 & 0.038 & 1.074 & 0.108 & 0.958 & 0.183 & 0.627 \\
      & 0.5 & $>0$ & 0.795 & 0.624 & 0.955 & 0.943 & 0.895 & 2.035 & 3.001 & 0.122 & 0.519 & 0.143 & 0.767 & 0.176 & 0.792 & 0.271 & 0.331  \\
      & 0.8 & $>0$ &  0.328 & 0.307 & 0.944 & 0.956 & 0.903 & 2.767 & 3.498 & 0.223 & 0.706 & 0.209 & 0.892 & 0.244 & 0.927 & 0.412 & 0.634  \\
      \hline 
      \multirow[t]{3}{*}{(500,550)} & 0.2 & $>0$ & 0.868 & 0.865 & 0.946 & 0.961 & 0.917 & 2.614 & 2.824 & 0.028 & 0.667 & 0.065 & 0.720 & 0.087 & 0.818 & 0.088 & 0.542 \\
      & 0.5 & $>0$ & 0.897 & 0.901 & 0.956 & 0.972 & 0.913 & 1.795 & 1.834 & 0.056 & 0.458 & 0.094 & 0.468 & 0.117 & 0.622 & 0.182 & 0.195 \\
      & 0.8 & $>0$ & 0.355 & 0.282 & 0.948 & 0.958 & 0.921 & 2.633 & 2.795 & -0.049 & 0.672 & 0.081 & 0.713 & 0.083 & 0.810 & 0.135 & 0.461 \\   
      \hline 
      \multirow[t]{3}{*}{(900,800)} & 0.2 & $>0$  & 0.971 & 0.963 & 0.952 & 0.958 & 0.923 & 1.988 & 2.395 & -0.015 & 0.507 & -0.021 & 0.611 & 0.053 & 0.773 & 0.073 & 0.314 \\
      & 0.5 & $>0$ & 0.984 & 0.981 & 0.953 & 0.953 & 0.928 & 1.383 & 1.380 & 0.046 & 0.353 & 0.058 & 0.352 & 0.072 & 0.525 & 0.082 & 0.142   \\
      & 0.8 & $>0$ & 0.407 & 0.385 & 0.946 & 0.952 & 0.922 & 1.991 & 2.082 & 0.062 & 0.508 & 0.059 & 0.531 & 0.069 & 0.614 & 0.095 & 0.344   \\
      \hline
    \end{tabular}
}
\end{table}

The simulation results are presented in Table \ref{tab1}. 
As the sample size increases, the proposed \texttt{IQTE} shows improved performance in both settings, evidenced by smaller biases and standard errors, and decreasing confidence interval lengths. This trend is further confirmed by the increasing power, measured by the rejection rates, for the settings where ${\Delta}_{\tau,\text {new}}>0$.
Compared to \texttt{Lasso}, \texttt{IQTE} exhibits smaller bias and larger variance in both settings. This bias-variance tradeoff aligns with our expectations. 
Additionally, \texttt{IQTE} generally has smaller bias and variance compared to both \texttt{RS} and \texttt{Deb}, particularly in the dense setting, thanks to the proposed variance-enhanced projection.
For coverage probabilities of confidence intervals, \texttt{IQTE} achieves the nominal 95\% confidence level in both settings, whereas the coverage probabilities for \texttt{Deb} are significantly lower than the nominal level. 
We also observe that the coverage probabilities for \texttt{RS} deviate from the nominal level in the dense loading vector setting, which is consistent with our theoretical findings.
Finally, in both settings, \texttt{IQTE} outperforms the most competitive method, \texttt{RS}, with smaller confidence interval lengths and higher power under the alternatives with $\Delta_{\tau,\text{new}}>0$.

\section{Real Data Analysis} \label{realdata.sec}
\subsection{Background} 
Numerous epidemiological studies have shown that serum lipid levels are crucial biomarkers for a variety of diseases, with low-density lipoprotein (LDL) and high-density lipoprotein (HDL) cholesterol serving as key clinical indicators. For instance, findings from the Framingham Heart Study indicate that both elevated LDL and reduced HDL cholesterol are associated with an increased risk of cardiovascular disease \citep{kannel1961}.

In recent years, a significant number of genome-wide association studies (GWAS) have investigated the genetic determinants of serum lipid levels \citep{thomas2018}. Beyond genetic factors, environmental influences, such as chemical contaminants and the climate, also play a substantial role in determining lipid levels, which prompts the emergence of environmental-wide association studies (EWAS) \citep{patel2010}. However, research examining the association between serum lipids and environmental factors remains comparatively limited. The National Health and Nutrition Examination Survey (NHANES) offers a valuable opportunity for researchers to address this need.
NHANES, conducted by the U.S. Centers for Disease Control and Prevention (CDC), monitors the health and nutritional status of the U.S. population, providing comprehensive data on various environmental factors. 
In a detailed analysis, \cite{patel2012} processes this data to identify the treatment effects of environmental volatile compounds on serum lipid levels. Among these volatile compounds, we focus on bromodichloromethane ($\text{CHBrCl}_2$), a byproduct of water chlorination linked to liver diseases \citep{zheng2021}. Our goal is to investigate the impact of $\text{CHBrCl}_2$ on HDL cholesterol as a biomarker for liver diseases \citep{fan2019}. 

\subsection{Data Analysis}
We conduct the analysis based on the research of \cite{patel2012}, utilizing laboratory and questionnaire data sourced from NHANES surveys conducted between 1999 and 2006. Following their initial selection of 127 covariates, we further refine the dataset by excluding variables with missing values, resulting in a final analysis comprising $p=101$ covariates. 
Examples include carotenoids, nutritional biochemistries and heavy metals; see a comprehensive list of these covariates in Section C 
of the Supplementary Material.

\begin{table}[t!]     
    \centering
    \caption{IQTE of $\text{CHBrCl}_2$ on HDL cholesterol level (g/dL) for five different individuals. Point estimators and 95\% confidence intervals are provided.\vspace{2mm}}
    \label{tab2}
    \resizebox{0.6\textwidth}{!}{
        \begin{tabular}{ccccc}
            \hline \textbf{SEQN}   &$\boldsymbol{\tau}$\quad &   \textbf{IQTE}    & \textbf{  Lower bound}   & \textbf{Upper bound} \vspace{1mm}\\
            \hline  
             \multirow[t]{3}{*}{30258} & 0.2 & 0.0036 & -0.0183&  0.0254  \\
             & 0.5 &   -0.0098 & -0.0219 &  0.0023  \\
             & 0.8 &   -0.0059 & -0.0267 &   0.0149 \\
             \hline 
              \multirow[t]{3}{*}{30519} & 0.2 & 0.0349  & -0.1087 & 0.1786  \\
              & 0.5 & 0.0100 & -0.1044 & 0.1244 \\	
              & 0.8 & 0.0646 &	-0.1969&  0.3260\\
              \hline
              \multirow[t]{3}{*}{30536} & 0.2 &  0.0196  & -0.0892&  0.1284 \\
              & 0.5 & -0.0042  & -0.1184 &  0.1099 \\		
              & 0.8 & -0.0253  & -0.1995&  0.1489 \\			
             \hline 
             \multirow[t]{3}{*}{30592} & 0.2 & 0.0060  &  -0.0240 & 0.0360 \\
             & 0.5 & 0.0071 & -0.0194  & 0.0337  \\		
             & 0.8 &0.0416 &   -0.0055 &  0.0886  \\			
            \hline
            \multirow[t]{3}{*}{30645} & 0.2 &  0.0093&  -0.0226  & 0.0413  \\
            & 0.5 & -0.0089  & -0.0333 & 0.0156 \\		
            & 0.8 & -0.0275  &-0.0738 &  0.0189 \\			
           \hline
            \end{tabular}
    }
\end{table}

The covariates are measured across $n=125$ subjects, comprising $n_1=49$ individuals exposed to the environment with high levels of $\text{CHBrCl}_2$ while $n_2=76$ individuals exposed to lower levels of this compound. 
Building upon the initial two samples, our objective is to investigate the IQTE of $\text{CHBrCl}_2$ on HDL cholesterol levels for five additional randomly selected subjects, identified by their respondent sequence numbers (SEQN) from NHANES: 30258, 30519, 30536, 30592, and 30645. Utilizing the method outlined in Algorithm \ref{alg1}, we derive the point estimates of IQTEs along with the corresponding 95\% confidence intervals at quantile levels $\tau=0.2,0.5\text{, and }0.8$. The results are presented in Table \ref{tab2}.
\begin{figure} 
    \caption{IQTE of $\text{CHBrCl}_2$ on HDL cholesterol (g/dL) for the  individual 30592. The 95\% confidence intervals are presented by the shadow area.} \label{fig1}
\centering
    \includegraphics[width=0.85\textwidth]{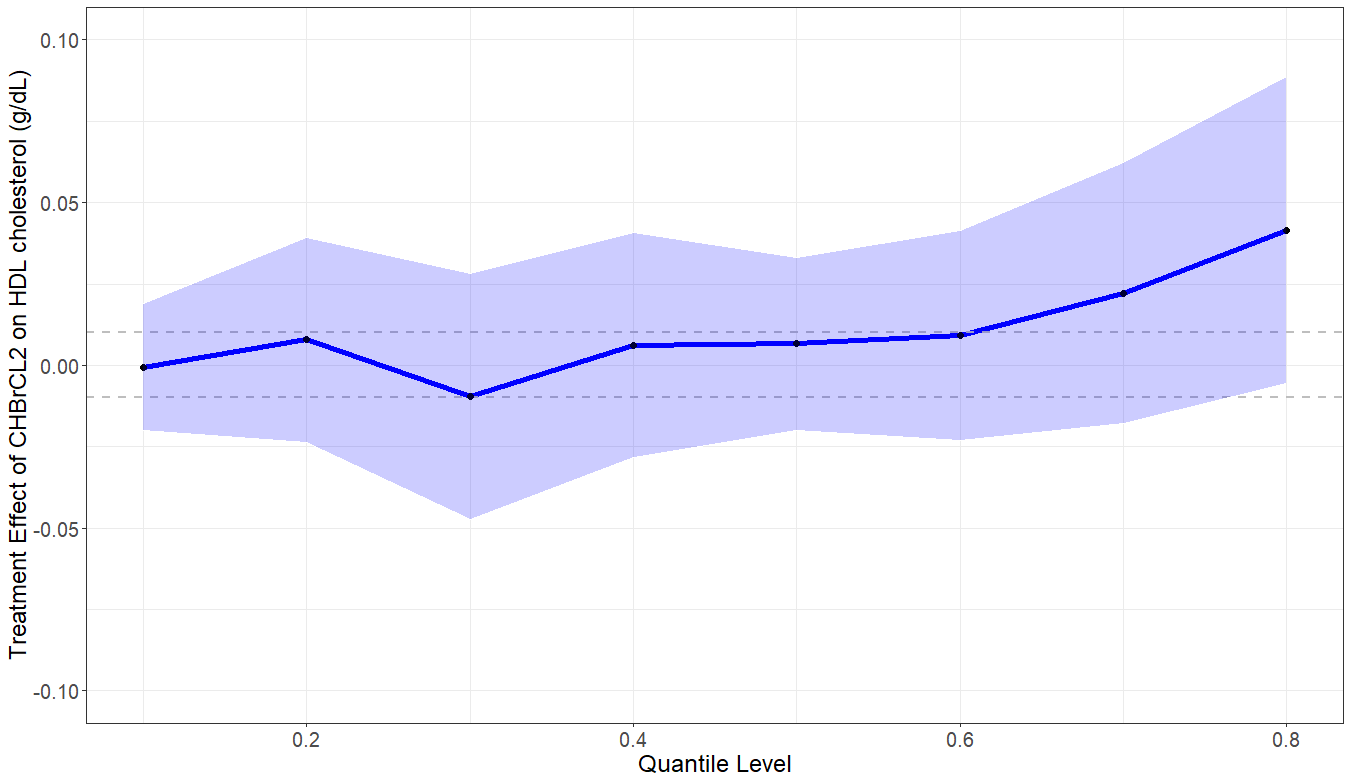}
\end{figure}

First, the estimates of treatment effects for the two individuals with SEQN 30519 and 30592 are positive at quantile levels $\tau=0.5$ and $0.8$, which contrast with those of the other three subjects. This reveals the variability in treatment effects across different individuals and affirms the existence of personalized heterogeneity. 
Second, except the individual 30258, the magnitude of treatment effects at $\tau=0.8$ significantly exceeds those at $\tau=0.2$ and $0.5$, which demonstrates substantial quantile heterogeneity. Figure \ref{fig1} visually depicts such heterogeneity by estimating the QTEs for individual 30592 at additional quantile levels, with the corresponding confidence intervals indicated by the shaded areas. It can be observed that for quantile levels $\tau \leq 0.6$, the treatment effects are relatively small, locating between $\pm 0.01$g/dL (the dashed lines), while the effects become more significant at upper quantiles. 
This numerical performance aligns with the biological theory of quantile-specific heritability of HDL cholesterol, which suggests that genetic markers become more detectable as HDL concentrations increase \citep[e.g.][]{williams}. In conclusion, the above observations underscore the importance of accounting for both personalized and quantile heterogeneities when analyzing treatment effects in practical applications.

\section{Discussion}\label{disc.sec}
In this article, we propose a debiased estimator for the IQTE that addresses both personalized and quantile heterogeneities. 
Leveraging a novel variance-enhancement constraint, this estimator and its corresponding inference procedures accommodate structure-free loading vectors. 
We validate the proposed procedures by examining the coverage probability of the confidence interval and the type I error in hypothesis testing. Additionally, we establish the minimax optimality of these procedures, demonstrating the superiority of our method.

For future research, exploring a convolution-type method, analogous to the approach in \cite{tan2022}, is promising. Such a method could simplify the computational challenges posed by the non-smoothness of the quantile loss function. However, the minimax optimality of inference procedures within this setting warrants further theoretical explorations.
In addition, our work can be extended to federated learning of IQTE across multiple sites in high-dimensional quantile regression settings. The development of integrative estimation and inference methods, along with the corresponding minimax results in this framework, remains an open area for future study.

\spacingset{1}
\normalem
\bibliographystyle{apalike}
\bibliography{refsqte}   

\begin{thebibliography}{}

\bibitem[Abadie et~al., 2002]{Abadie2002}
Abadie, A., Angrist, J., and Imbens, G. (2002).
\newblock Instrumental variables estimates of the effect of subsidized training
  on the quantiles of trainee earnings.
\newblock {\em Econometrica}, 70(1):91--117.

\bibitem[Abadie and Imbens, 2006]{Abadie2006}
Abadie, A. and Imbens, G.~W. (2006).
\newblock Large sample properties of matching estimators for average treatment
  effects.
\newblock {\em Econometrica}, 74(1):235--267.

\bibitem[Angrist et~al., 2006]{angrist2006}
Angrist, J., Chernozhukov, V., and Fernández-Val, I. (2006).
\newblock Quantile regression under misspecification, with an application to
  the {U.S.} wage structure.
\newblock {\em Econometrica}, 74(2):539--563.

\bibitem[Athey et~al., 2018]{athey2018}
Athey, S., Imbens, G.~W., and Wager, S. (2018).
\newblock Approximate residual balancing: debiased inference of average
  treatment effects in high dimensions.
\newblock {\em J. R. Stat. Soc. Ser. B. Stat. Methodol.}, 80(4):597--623.

\bibitem[Belloni and Chernozhukov, 2011]{Belloni2011}
Belloni, A. and Chernozhukov, V. (2011).
\newblock {{$\ell_1$}-penalized quantile regression in high-dimensional sparse
  models}.
\newblock {\em Ann. Statist.}, 39(1):82--130.

\bibitem[Belloni et~al., 2019]{belloni2019}
Belloni, A., Chernozhukov, V., and Kato, K. (2019).
\newblock Valid post-selection inference in high-dimensional approximately
  sparse quantile regression models.
\newblock {\em J. Amer. Statist. Assoc.}, 114(526):749--758.

\bibitem[Beyerlein et~al., 2011]{beyerlein2011}
Beyerlein, A., von Kries~R, Ness, A., and Ong, K. (2011).
\newblock Genetic markers of obesity risk: stronger associations with body
  composition in overweight compared to normal-weight children.
\newblock {\em PLoS One}, 6(4):e19057.

\bibitem[Bickel et~al., 2009]{Bickel2009}
Bickel, P.~J., Ritov, Y., and Tsybakov, A.~B. (2009).
\newblock Simultaneous analysis of lasso and {Dantzig} selector.
\newblock {\em Ann. Statist.}, 37(4):1705--1732.

\bibitem[Bradic et~al., 2011]{bradic2011}
Bradic, J., Fan, J., and Wang, W. (2011).
\newblock Penalized composite quasi-likelihood for ultrahigh dimensional
  variable selection.
\newblock {\em J. R. Stat. Soc. Ser. B. Stat. Methodol.}, 73(3):325–349.

\bibitem[Bradic et~al., 2024]{bradic2023highdimensional}
Bradic, J., Ji, W., and Zhang, Y. (2024).
\newblock {High-dimensional inference for dynamic treatment effects}.
\newblock {\em Ann. Statist.}, 52(2):415--440.

\bibitem[Bradic and Kolar, 2017]{bradic2017uniform}
Bradic, J. and Kolar, M. (2017).
\newblock Uniform inference for high-dimensional quantile regression: linear
  functionals and regression rank scores.
\newblock {\em arXiv:1702.06209}.

\bibitem[Briollais and Durrieu, 2014]{briollais2014}
Briollais, L. and Durrieu, G. (2014).
\newblock Application of quantile regression to recent genetic and-omic
  studies.
\newblock {\em Human genetics}, 133(8):951--966.

\bibitem[Bronnenberg et~al., 2012]{Bart2012}
Bronnenberg, B.~J., Dubé, J.-P.~H., and Gentzkow, M. (2012).
\newblock The evolution of brand preferences: {evidence} from consumer
  migration.
\newblock {\em Am. Econ. Rev.}, 102(6):2472--2508.

\bibitem[B{\"u}hlmann and van~de Geer, 2011]{peter2011}
B{\"u}hlmann, P. and van~de Geer, S. (2011).
\newblock {\em Statistics for high-dimensional data: methods, theory and
  applications}.
\newblock Springer Science \& Business Media.

\bibitem[Cai et~al., 2021]{cai2021}
Cai, T., Cai, T.~T., and Guo, Z. (2021).
\newblock Optimal statistical inference for individualized treatment effects in
  high-dimensional models.
\newblock {\em J. R. Stat. Soc. Ser. B. Stat. Methodol.}, 83(4):669--719.

\bibitem[Cai et~al., 2011]{cai2011}
Cai, T., Tian, L., Wong, P.~H., and Wei, L. (2011).
\newblock Analysis of randomized comparative clinical trial data for
  personalized treatment selections.
\newblock {\em Biostatistics}, 12(2):270--282.

\bibitem[Cai and Guo, 2017]{cai2017}
Cai, T.~T. and Guo, Z. (2017).
\newblock {Confidence intervals for high-dimensional linear regression: minimax
  rates and adaptivity}.
\newblock {\em Ann. Statist.}, 45(2):615--646.

\bibitem[Cai et~al., 2023]{cai2023statistical}
Cai, T.~T., Guo, Z., and Xia, Y. (2023).
\newblock Statistical inference and large-scale multiple testing for
  high-dimensional regression models.
\newblock {\em TEST}, 32:1135--1171.

\bibitem[Cai and Zhou, 2012]{Cai_2012}
Cai, T.~T. and Zhou, H.~H. (2012).
\newblock Optimal rates of convergence for sparse covariance matrix estimation.
\newblock {\em Ann. Statist.}, 40(5):2389--2420.

\bibitem[Candès and Tao, 2007]{candes2007}
Candès, E. and Tao, T. (2007).
\newblock {The Dantzig selector: statistical estimation when {$p$} is much
  larger than {$n$}}.
\newblock {\em Ann. Statist.}, 35(6):2313--2351.

\bibitem[Chao et~al., 2017]{chao2017}
Chao, S.-K., Volgushev, S., and Cheng, G. (2017).
\newblock Quantile processes for semi and nonparametric regression.
\newblock {\em Electron. J. Stat.}, 11(2):3272--3331.

\bibitem[Cheng et~al., 2022]{chao2022}
Cheng, C., Feng, X., Huang, J., and Liu, X. (2022).
\newblock Regularized projection score estimation of treatment effects in
  high-dimensional quantile regression.
\newblock {\em Statist. Sinica}, 32:23--41.

\bibitem[Chernozhukov and Hansen, 2005]{Chernozhukov2005}
Chernozhukov, V. and Hansen, C. (2005).
\newblock An {IV} model of quantile treatment effects.
\newblock {\em Econometrica}, 73(1):245--261.

\bibitem[Fan et~al., 2014]{fan2014}
Fan, J., Fan, Y., and Barut, E. (2014).
\newblock Adaptive robust variable selection.
\newblock {\em Ann. Statist.}, 42(1):324–351.

\bibitem[Fan et~al., 2019]{fan2019}
Fan, N., Peng, L., Xia, Z., Zhang, L., Song, Z., Wang, Y., and Peng, Y. (2019).
\newblock {Triglycerides to high-density lipoprotein cholesterol ratio as a
  surrogate for nonalcoholic fatty liver disease: a cross-sectional study}.
\newblock {\em Lipids Health Dis.}, 18(39):1--6.

\bibitem[Fan and Liu, 2016]{fan2016}
Fan, Y. and Liu, R. (2016).
\newblock A direct approach to inference in nonparametric and semiparametric
  quantile models.
\newblock {\em J. Econom.}, 191(1):196--216.

\bibitem[Giessing and Wang, 2023]{giessing2023debiased}
Giessing, A. and Wang, J. (2023).
\newblock Debiased inference on heterogeneous quantile treatment effects with
  regression rank-scores.
\newblock {\em J. R. Stat. Soc. Ser. B. Stat. Methodol.}, 85(5):1561--1588.

\bibitem[Goldhirsch et~al., 2002]{gold2002}
Goldhirsch, A., Coates, A., Collins, J., Thurlimann, B., Senn, H.-J., Holmberg,
  S., Lindtner, J., Veronesi, A., Cortés-Funes, H., Castiglione-Gertsch, M.,
  Nasi, M., Egli, G., Rabaglio, M., Maibach, R., Gerber, M., Hiltbrunner, A.,
  Gelber, R., Price, K., Bonetti, M., and Ferrero, A. (2002).
\newblock Endocrine responsiveness and tailoring adjuvant therapy for
  postmenopausal lymph node-negative breast cancer: {a} randomized trial.
\newblock {\em J. Nat. Cancer Inst.}, 94(14):1054--1065.

\bibitem[Gutenbrunner and Jureckova, 1992]{guten1992}
Gutenbrunner, C. and Jureckova, J. (1992).
\newblock Regression rank scores and regression quantiles.
\newblock {\em Ann. Statist.}, 20(1):305--330.

\bibitem[He et~al., 2023]{he2023}
He, X., Pan, X., Tan, K.~M., and Zhou, W. (2023).
\newblock Smoothed quantile regression with large-scale inference.
\newblock {\em J. Econom.}, 232(2):367--388.

\bibitem[Hirano et~al., 2003]{Hirano2003}
Hirano, K., Imbens, G., and Ridder, G. (2003).
\newblock Efficient estimation of average treatment effects using the estimated
  propensity score.
\newblock {\em Econometrica}, 71(4):1161--1189.

\bibitem[Hoffmann et~al., 2018]{thomas2018}
Hoffmann, T.~J., Theusch, E., Haldar, T., Ranatunga, D.~K., Jorgenson, E.,
  Medina, M.~W., Kvale, M.~N., Kwok, P.-Y., Schaefer, C., Krauss, R.~M.,
  Iribarren, C., and Risch, N. (2018).
\newblock A large electronic-health-record-based genome-wide study of serum
  lipids.
\newblock {\em Nature Genetics}, 50:401--413.

\bibitem[Imai and Ratkovic, 2013]{kosuke2013}
Imai, K. and Ratkovic, M. (2013).
\newblock Estimating treatment effect heterogeneity in randomized program
  evaluation.
\newblock {\em Ann. Appl. Stat.}, 7(1):443--470.

\bibitem[Imbens and Wooldridge, 2009]{imbens2009}
Imbens, G.~W. and Wooldridge, J.~M. (2009).
\newblock Recent developments in the econometrics of program evaluation.
\newblock {\em J. Econ. Lit.}, 47(1):5--86.

\bibitem[Javanmard and Montanari, 2014]{JMLR:v15:javanmard14a}
Javanmard, A. and Montanari, A. (2014).
\newblock Confidence intervals and hypothesis testing for high-dimensional
  regression.
\newblock {\em J. Mach. Learn. Res.}, 15(82):2869--2909.

\bibitem[Kannel et~al., 1961]{kannel1961}
Kannel, W.~B., Dawber, T.~R., Kagan, A., Revotskie, N., and Stokes, J. (1961).
\newblock Factors of risk in the development of coronary heart disease—six
  year follow-up experience: the {Framingham} study.
\newblock {\em Ann. Intern. Med.}, 55:33--50.

\bibitem[Koenker, 2005]{koenker_2005}
Koenker, R. (2005).
\newblock {\em Quantile regression}.
\newblock Cambridge University Press.

\bibitem[Koenker, 2010]{koenker2010}
Koenker, R. (2010).
\newblock Rank tests for heterogeneous treatment effects with covariates.
\newblock {\em Inst. Math. Stat. Collect.}, 7:134--142.

\bibitem[Koenker, 2017]{koenker2017}
Koenker, R. (2017).
\newblock Quantile regression: 40 years on.
\newblock {\em Annu. Rev. Econ.}, 9(1):155--176.

\bibitem[Koenker and Xiao, 2002]{koenker2002process}
Koenker, R. and Xiao, Z. (2002).
\newblock Inference on the quantile regression process.
\newblock {\em Econometrica}, 70(4):1583--1612.

\bibitem[Kosorok and Laber, 2019]{kosorok2019}
Kosorok, M.~R. and Laber, E.~B. (2019).
\newblock Precision medicine.
\newblock {\em Annu. Rev. Stat. Appl.}, 6(1):263--286.

\bibitem[Liu et~al., 2017]{liu2017}
Liu, W., Zhang, Z., Nie, L., and Soon, G. (2017).
\newblock A case study in personalized medicine: {rilpivirine} versus efavirenz
  for treatment-naive {HIV} patients.
\newblock {\em J. Amer. Statist. Assoc.}, 112(520):1381--1392.

\bibitem[McCaw et~al., 2020]{mccaw2020}
McCaw, Z.~R., Tian, L., Vassy, J.~L., Ritchie, C.~S., Lee, C.~C., Kim, D.~H.,
  and Wei, L.~J. (2020).
\newblock How to quantify and interpret treatment effects in comparative
  clinical studies of {COVID}-19.
\newblock {\em Ann. Intern. Med.}, 173(8):632--637.

\bibitem[Patel et~al., 2010]{patel2010}
Patel, C.~J., Bhattacharya, J., and Butte, A.~J. (2010).
\newblock An environment-wide association study ({EWAS}) on type 2 diabetes
  mellitus.
\newblock {\em PloS one}, 5(5):e10746.

\bibitem[Patel et~al., 2012]{patel2012}
Patel, C.~J., Cullen, M.~R., Ioannidis, J.~P., and Butte, A.~J. (2012).
\newblock {Systematic evaluation of environmental factors: persistent
  pollutants and nutrients correlated with serum lipid levels}.
\newblock {\em Int. J. Epidemiol.}, 41(3):828–843.

\bibitem[Sabine, 2005]{sabine2005}
Sabine, C. (2005).
\newblock {AIDS} events among individuals initiating {HAART}: Do some patients
  experience a greater benefit from {HAART} than others?
\newblock {\em AIDS}, 19(17):1995--2000.

\bibitem[Sasaki and Ura, 2023]{SASAKI2023394}
Sasaki, Y. and Ura, T. (2023).
\newblock Estimation and inference for policy relevant treatment effects.
\newblock {\em J. Econom.}, 234(2):394--450.

\bibitem[Sergio, 2007]{firpo2007}
Sergio, F. (2007).
\newblock Efficient semiparametric estimation of quantile treatment effects.
\newblock {\em Econometrica}, 75(1):259--276.

\bibitem[Sun and He, 2021]{Sun2021}
Sun, Y. and He, X. (2021).
\newblock Model-based bootstrap for detection of regional quantile treatment
  effects.
\newblock {\em J. Nonparametric Stat.}, 33(2):299--320.

\bibitem[Tan et~al., 2021]{tan2022}
Tan, K.~M., Wang, L., and Zhou, W.-X. (2021).
\newblock High-dimensional quantile regression: convolution smoothing and
  concave regularization.
\newblock {\em J. R. Stat. Soc. Ser. B. Stat. Methodol.}, 84(1):205--233.

\bibitem[Tibshirani, 1996]{lasso}
Tibshirani, R. (1996).
\newblock Regression shrinkage and selection via the lasso.
\newblock {\em J. R. Stat. Soc. Ser. B. Stat. Methodol.}, 58(1):267--288.

\bibitem[Vaart and Wellner, 1996]{vaart1996}
Vaart, A.~W. and Wellner, J.~A. (1996).
\newblock {\em Weak convergence and empirical processes}.
\newblock Springer Series in Statistics.

\bibitem[van~de Geer et~al., 2014]{vandegeer2014}
van~de Geer, S., B{\"u}hlmann, P., Ritov, Y., and Dezeure, R. (2014).
\newblock {On asymptotically optimal confidence regions and tests for
  high-dimensional models}.
\newblock {\em Ann. Statist.}, 42(3):1166--1202.

\bibitem[Volgushev et~al., 2019]{chao2019}
Volgushev, S., Chao, S.-K., and Cheng, G. (2019).
\newblock Distributed inference for quantile regression process.
\newblock {\em Ann. Statist.}, 47(3):1634--1662.

\bibitem[Wang et~al., 2009]{wang2009}
Wang, H.~J., Zhu, Z., and Zhou, J. (2009).
\newblock Quantile regression in partially linear varying coefficient models.
\newblock {\em Ann. Statist.}, 37(6B):3841--3866.

\bibitem[Wang et~al., 2020]{wang2020}
Wang, J., He, X., and Xu, G. (2020).
\newblock Debiased inference on treatment effect in a high-dimensional model.
\newblock {\em J. Amer. Statist. Assoc.}, 115(529):442--454.

\bibitem[Wang et~al., 2018a]{lanwang2018}
Wang, L., Zhou, Y., Song, R., and Sherwood, B. (2018a).
\newblock Quantile-optimal treatment regimes.
\newblock {\em J. Amer. Statist. Assoc.}, 113(523):1243--1254.

\bibitem[Wang et~al., 2022]{wang2022}
Wang, T., Ionita-Laza, I., and Wei, Y. (2022).
\newblock Integrated quantile rank test {(iQRAT)} for gene-level associations.
\newblock {\em Ann. Appl. Stat.}, 16(3):1423--1444.

\bibitem[Wang et~al., 2018b]{wang2018}
Wang, Y., Fu, H., and Zeng, D. (2018b).
\newblock Learning optimal personalized treatment rules in consideration of
  benefit and risk: {with} an application to treating type 2 diabetes patients
  with insulin therapies.
\newblock {\em J. Amer. Statist. Assoc.}, 113(521):1--13.

\bibitem[Williams, 2020]{williams}
Williams, P.~T. (2020).
\newblock Quantile-specific heritability of high-density lipoproteins with
  implications for precision medicine.
\newblock {\em J. Clin. Lipidol.}, 14(4):448--458.

\bibitem[Wu et~al., 2023]{wu2023}
Wu, Y., Wang, L., and Fu, H. (2023).
\newblock Model-assisted uniformly honest inference for optimal treatment
  regimes in high dimension.
\newblock {\em J. Amer. Statist. Assoc.}, 118(541):305--314.

\bibitem[Xia et~al., 2018]{xia2018}
Xia, Y., Cai, T.~T., and Li, H. (2018).
\newblock {Joint testing and false discovery rate control in high-dimensional
  multivariate regression}.
\newblock {\em Biometrika}, 105(2):249--269.

\bibitem[Yan et~al., 2023]{yan2023}
Yan, Y., Wang, X., and Zhang, R. (2023).
\newblock Confidence intervals and hypothesis testing for high-dimensional
  quantile regression: convolution smoothing and debiasing.
\newblock {\em J. Mach. Learn. Res.}, 24(245):1--49.

\bibitem[Yang and He, 2012]{yang2012}
Yang, Y. and He, X. (2012).
\newblock {Bayesian empirical likelihood for quantile regression}.
\newblock {\em Ann. Statist.}, 40(2):1102--1131.

\bibitem[Yang et~al., 2016]{yang2016}
Yang, Y., Wang, H.~J., and He, X. (2016).
\newblock Posterior inference in {Bayesian} quantile regression with asymmetric
  {Laplace} likelihood.
\newblock {\em Int. Statist. Rev.}, 84(3):327--344.

\bibitem[Zhang and Zhang, 2014]{zhang2014}
Zhang, C.~H. and Zhang, S.~S. (2014).
\newblock Confidence intervals for low dimensional parameters in high
  dimensional linear models.
\newblock {\em J. R. Stat. Soc. Ser. B. Stat. Methodol.}, 76(1):217--242.

\bibitem[Zhang et~al., 2020]{zhang2020}
Zhang, Y., Wang, L., Yu, M., and Shao, J. (2020).
\newblock Quantile treatment effect estimation with dimension reduction.
\newblock {\em Stat. Theory Relat. Fields.}, 4(2):202--213.

\bibitem[Zhao et~al., 2015]{zhao2015general}
Zhao, T., Kolar, M., and Liu, H. (2015).
\newblock A general framework for robust testing and confidence regions in
  high-dimensional quantile regression.
\newblock {\em arXiv:1412.8724}.

\bibitem[Zhao et~al., 2020]{zhao2020}
Zhao, W., Zhang, F., and Lian, H. (2020).
\newblock Debiasing and distributed estimation for high-dimensional quantile
  regression.
\newblock {\em IEEE Trans. Neural Netw. Learn. Syst.}, 31(7):2569--2577.

\bibitem[Zheng et~al., 2015]{zheng2015}
Zheng, Q., Peng, L., and He, X. (2015).
\newblock {Globally adaptive quantile regression with ultra-high dimensional
  data}.
\newblock {\em Ann. Statist.}, 43(5):2225--2258.

\bibitem[Zheng et~al., 2021]{zheng2021}
Zheng, S., Yang, Y., Wen, C., Liu, W., Cao, L., Feng, X., Chen, J., Wang, H.,
  Tang, Y., Tian, L., Wang, X., and Yang, F. (2021).
\newblock Effects of environmental contaminants in water resources on
  nonalcoholic fatty liver disease.
\newblock {\em Environ. Int.}, 154:106555.

\end{thebibliography}

\end{document}